# ASYMPTOTICS OF THE PROBABILITY MINIMIZING A "DOWN-SIDE" RISK

By Hiroaki Hata, Hideo Nagai[1] and Shuenn-Jyi Sheu[2]

*Academia Sinica, Osaka University and Academia Sinica*

We consider a long-term optimal investment problem where an investor tries to minimize the probability of falling below a target growth rate. From a mathematical viewpoint, this is a large deviation control problem. This problem will be shown to relate to a risk-sensitive stochastic control problem for a sufficiently large time horizon. Indeed, in our theorem we state a duality in the relation between the above two problems. Furthermore, under a multidimensional linear Gaussian model we obtain explicit solutions for the primal problem.

**1. Introduction.** In recent studies of finance, it has been of great concern to consider problems from risk management. In this paper, we consider the problem of the control of a down-side risk probability for an investor. To minimize such probabilities and obtain an optimal (or nearly optimal) portfolio, we relate the problem of the portfolio optimization with risk sensitive criterion. In [35, 36], from the consideration of a performance index for funds when compared to a benchmark, a similar problem is considered. In [30, 31], problems of maximizing the up-side chance probability are studied. These studies show potential applications of risk-sensitive dynamic management to the problem of risk management.

Portfolio optimization with risk-sensitive criterion has been considered in several recent works (see [2–4, 7, 10, 14–16, 18, 19, 21, 25, 28] and [29]). The problem is to maximize the expected utility of terminal wealth with the HARA utility function being considered,

$$(1.1) \qquad \max\left\{ E\left[\frac{1}{\gamma}(X_T^\pi)^\gamma\right]\right\},$$

Received August 2007; revised September 2008.

[1]Supported in part by Grant-in-Aid for scientific research 20340019 JSPS.

[2]Supported by NSC Grant 96-2119-M-001-002.

where $X_t^\pi$ is the wealth process using strategy $\pi$, and the maximum is taken over a class of admissible strategies. In the following discussion, it is also convenient to replace (1.1) by

$$(1.1a) \qquad \max\{E[(X_T^\pi)^\gamma]\}, \qquad \gamma > 0,$$

and

$$(1.1b) \qquad \min\{E[(X_T^\pi)^\gamma]\}, \qquad \gamma < 0.$$

The parameter $\gamma$ is taken from $(-\infty, 1)$ and for $\gamma = 0$, $\frac{1}{\gamma}(X_T^\pi)^\gamma$ is interpreted as $\log X_T^\pi$ which is a Kelly utility. When the HARA utility is used, the problem can be reformulated as a risk-sensitive control problem where $\gamma$ plays the role of risk-sensitive parameter. Then methods from the theory of stochastic control can be applied to the portfolio optimization problem. In particular, we may apply dynamic programming to solve the problem. This approach is very useful in discrete time models (see [24] for some idea). In continuous time, we need to solve the Bellman equation. In the past few years, a special class of models has been extensively studied in this frame work where a factor process is introduced and the return, together with the volatility of stock prices, is affected by the factor process. That is, we assume in a market that we have one bank account with price $S_t^0$, $m$ stocks with prices $S_t^i, i = 1, \ldots, m$, and $n$ economic factors $Y_t = (Y_t^1, Y_t^2, \ldots, Y_t^n)$. Their dynamics are given by

$$(1.2) \qquad dS_t^0 = r(Y_t)S_t^0\,dt, \qquad S_0^0 = s^0,$$

$$(1.3) \qquad dS_t^i = S_t^i\left\{\alpha^i(Y_t)\,dt + \sum_{k=1}^{n+m} \sigma_k^i(Y_t)\,dW_t^k\right\}, \qquad S^i(0) = s^i, i = 1, \ldots, m,$$

$$(1.4) \qquad dY_t = \beta(Y_t)\,dt + \lambda(Y_t)\,dW_t, \qquad Y(0) = y \in \mathbb{R}^n,$$

where $W_t = (W_t^k)_{k=1,\ldots,(n+m)}$ is an $(n+m)$-dimensional standard Brownian motion defined on a filtered probability space $(\Omega, \mathcal{F}, P, \mathcal{F}_t)$. A closed form optimal strategy may be obtained by solving the Bellman equation. See [10] for the initial study and the subsequence works cited above. For an introduction of the theory of the risk-sensitive control problem, one can see [37]. For the mathematical theory of the risk-sensitive control problem and the connection with the robust control problem, one can see [11, 12, 17, 27].

Bielecki and Pliska ([2], Section 6) mention a possible use of portfolio optimization with risk-sensitive criterion in the study of the problem of upside chance and down-side risk. The problem of up-side chance is to consider

$$(1.5) \qquad \max\left\{P\left(\frac{1}{T}\log X_T^\pi \geq c\right)\right\}$$

for large $T$ where the maximization is taken for $\pi$ in a class of admissible strategies. Or more generally,

$$(1.6) \qquad \max\left\{P\left(\frac{1}{T}\log\frac{X_T^\pi}{I_T}\geq c\right)\right\},$$

where $I_T$ is a benchmark process. For simplicity, we take $I_t = 1$ in the following discussion. A mathematical theory was developed in [30] and [31] for the maximization of the up-side chance probabilities. It is shown that the probabilities in (1.5) are related to those in (1.1) with $0 < \gamma < 1$ through a duality relation. A nearly optimal portfolio for (1.5) can be obtained from an optimal portfolio for (1.1) for a particular chosen $\gamma$. Some idea from large deviation theory is implicitly used in his approach. This will be described later in this section. Interestingly, the maximization of the probability of up-side chance is not a conventional optimization problem and, in general, is difficult to treat. Therefore, studies in [30, 31] suggest the possibility of indirectly using the theory of stochastic control in such a nonconventional optimization problem. See also [19] for more studies on this problem. However, Sekine [33] recently tried to use a duality approach to treat such problems. In [30], Pham also proposed to develop a mathematical theory for the down-side risk probability,

$$(1.7) \qquad \min\left\{P\left(\frac{1}{T}\log X_T^\pi \leq k\right)\right\},$$

or more generally,

$$(1.8) \qquad \min\left\{P\left(\frac{1}{T}\log\frac{X_T^\pi}{I_T}\leq k\right)\right\}.$$

The problem is to minimize (1.6) and obtain an optimal (or nearly optimal) portfolio. Here $k$ is considered such that the event has small probabilities. Hence we are dealing with a rare event. From the consideration of finance, such rare events in down-side risk are not favorable to an investor. Therefore its occurrence may result in a significant consequence in portfolio management. Hence the study of (1.7) or (1.8) may be of meaningful implications in finance. See interesting discussions in [5, 6] and [35, 36] for some related problems where the consideration is the performance index for funds. In this paper, we develop some mathematical analyses for (1.7). Similar to [30, 31], we will show a dual relation between (1.1) and (1.7) for large $T$. The result says that for $k$, there is a correspondence $\gamma(k) < 0$ such that an optimal portfolio of (1.1) with $\gamma = \gamma(k)$ is a nearly optimal portfolio of (1.7). The meaning of this result is that an investor who wants to control (1.7) for a particular $k$ will have the same behavior as an investor whose risk parameter is $\gamma(k)$. See a similar result in [36] for some discrete time models.

In this paper we consider only the linear Gaussian models. But our method may be applied to more general (nonlinear) models. Gaussian models have practical uses for practitioners. A simple Gaussian factor model was first proposed in [26]. Such models have several important properties. It is much easier to estimate the coefficients by using linear regression (see [32]). This has practical applications. Tractibility is another big advantage of these models. For such models, a closed-form solution for the optimal portfolio selection problem (1.1) is possible by solving a Ricatti equation, which is a matrix equation and is easier to solve [[2, 3], [14–16, 25]]. For Gaussian models, we have some finer mathematical results [25]. In the later sections, we will see that these results will be crucial for our analysis. It is important to consider the model when some factors cannot be observed directly. The solution is far away from complete. One can find some results for linear Gaussian models in [18, 29, 34]. However, this is not our main concern in this paper.

The papers [30, 31] (and also this paper) consider (1.5) [and (1.7), respectively] for such $c$ (and $k$) that (1.5) [and (1.7)] has small probability. That is, we are dealing with large deviation probabilities. Therefore, we expect that some idea from large deviation theory can be used to relate (1.1) and (1.5) [or (1.7)] that we now explain. The formal calculation given in the following may be instructive to see the idea. For a given strategy $\pi$, assume $Z_T = \log X_T^\pi$ satisfies a large deviation principle with rate $I(k, \pi)$. Formally, this means

$$P\left(\frac{\log X_T^\pi}{T} \simeq k\right) \simeq \exp(I(k, \pi)T)$$

as $T \to \infty$. The Laplace–Varadhan lemma (see [8]) implies

$$E[(X_T^\pi)^\gamma] = E[\exp(\gamma \log X_T^\pi)] \simeq \exp(T\Phi(\gamma, \pi)),$$

where

$$\Phi(\gamma, \pi) = \sup_k \{\gamma k + I(k, \pi)\}.$$

If we want to minimize $\Phi(\gamma, \pi)$ [and minimize $I(k, \pi)$] over $\pi$, then

(1.9)  $$\inf_\pi \Phi(\gamma, \pi) = \inf_\pi \sup_k \{\gamma k + I(k, \pi)\}.$$

Denote

$$J(k) = \inf_\pi \{I(k, \pi)\},$$

$$\Phi(\gamma) = \inf_\pi \{\Phi(\gamma, \pi)\}.$$

If we are allowed to change the order of inf and sup on the right-hand side of (1.9), then we obtain

$$(1.10) \qquad \Phi(\gamma) = \sup_k \{\gamma k + J(k)\}.$$

If $J(k)$ is a concave function, then we expect to have the dual relation

$$(1.11) \qquad J(k) = \inf_\gamma \{\Phi(\gamma) - \gamma k\}$$

(see [9]). On the other hand, if we want to maximize $\Phi(\gamma, \pi)$ [and maximize $I(k, \pi)$] over possible $\pi$, then

$$(1.12) \qquad \sup_\pi \Phi(\gamma, \pi) = \sup_\pi \sup_k \{\gamma k + I(k, \pi)\}.$$

Denote

$$I(k) = \sup_\pi \{I(k, \pi)\},$$

$$\overline{\Phi}(\gamma) = \sup_\pi \{\Phi(\gamma, \pi)\}.$$

Then

$$\overline{\Phi}(\gamma) = \sup_k \{\gamma k + I(k)\}.$$

If $I(k)$ is concave function, then we expect to have the dual relation,

$$I(k) = \inf_\gamma \{\overline{\Phi}(\gamma) - \gamma k\}.$$

In [30, 31], through the above intuition, the following relation was proved:

$$\overline{\Pi}(c) = \inf_{0 < \gamma < 1} \{\Psi(\gamma) - \gamma c\},$$

where, for large $T$,

$$\text{value of (1.1)} \simeq \frac{1}{\gamma} \exp(T\Psi(\gamma)),$$

$$\text{value of (1.5)} \simeq \exp(T\overline{\Pi}(c)).$$

On the other hand, one of the main results of our paper is to prove

$$\Pi(k) = \inf_{\gamma < 0} \{\Psi(\gamma) - \gamma k\},$$

where, for large $T$,

$$\text{value of (1.7)} \simeq \exp(T\Pi(k)).$$

The problem analyzed here is closely related to the problem studied in [30, 31], since

$$P\left(\frac{1}{T}\log X_T^\pi \le k\right) = 1 - P\left(\frac{1}{T}\log X_T^\pi > k\right).$$

However, in both studies we are interested in the region of $c$ and $k$ such that probabilities in (1.5) and (1.7) are small. This means in this paper we study the region of $c$ such that (1.5) has large probability while in [30, 31], it is studied in the region of $c$ that (1.5) has small probability. This explains why the results of [30, 31] cannot be readily applied to the calculation of the down-side risk probabilities. We consider both of the problems for large $T$. We show in this paper that the minimization of (1.7) relates to the problem (1.1) for $\gamma < 0$. However, the probability in (1.5) relates to problem (1.1) for $0 < \gamma < 1$ as was shown in [30]. We will show that an optimal (or nearly optimal) portfolio for (1.7) can be derived from an optimal portfolio of (1.1) for a proper chosen $\gamma < 0$ that relates to $k$ through the duality relation. This is expected, after [30, 31]. However, our result does not follow from a simple Markov inequality as in [30, 31]. We need to use some finer results for (1.1), taken from [25] and some extensions given in Section 2 later. This is unexpected from the argument in [30, 31] and is also unexpected from the large deviation theory. Intuitively, when changing the order of inf and sup in (1.9) to obtain (1.10), it may cause some difficulty. Mathematically, we also want to mention another important point. The convexity of the value of (1.1b) with respect to $\gamma$ after taking limit for large $T$ will be crucial in our analysis. We note that the convexity of (1.1b) does not follow for finite $T$. This property is easily seen for (1.1a) from our formal derivation given in the previous paragraph, since taking the maximum of a family of convex functions is a convex function.

The paper is organized as follows. In Section 2, the model will be explicitly stated and the mathematical problem of down-side risk will be formulated. The associated problems of portfolio optimization will also be described. Since the problem has different formulations in an infinite time horizon, we will careful state our problem. We also give some remarks to mention other related problems. The main results will be stated in this section. These include the main theorem (Theorem 2.1), and several important propositions (Propositions 2.1–2.7). The results are presented in a way that a proof of Theorem 2.1 will follow using mainly results in the propositions and some minor properties. The proof of our main theorem also suggests that our main result will follow when we can prove the statements given in several crucial propositions. Therefore, the approach in this paper can be applied to more general models (including some nonlinear models). The problem for some nonlinear models is currently under investigation.

In Section 3, we give the proof of our main theorem. In Section 4, we prove the propositions stated in Section 2. These propositions concern the portfolio optimization problem (1.1) and the behaviors of the solution of the Bellman equation. Some results in these propositions have been obtained in [25]. They include Propositions 2.1, 2.2 and 2.3. Therefore, in the proof of these results, we give a sketch or modification of the arguments in [25] and mention the places in [25] where one can find the details of analysis used. Propositions 2.4, 2.5, 2.6 and 2.7 are some new results that are suggested by the problem studied in this paper. Therefore, some more details of the proofs will be given for these propositions.

## 2. The problem and main results.

2.1. *Down-side risk problem.* We consider a market model consisting of one bank account and $m$ risky stocks. The interest rate of bank account as well as the returns and volatility of stocks are affected by $n$ economic factors. The price of bank account $S^0$ and the price of risky stocks $S^i, i = 1, \ldots, m$, are given by

$$(2.1) \quad dS_t^0 = r(Y_t)S_t^0 \, dt, \qquad S_0^0 = s^0,$$

$$(2.2) \quad dS_t^i = S_t^i \left\{ \alpha^i(Y_t) \, dt + \sum_{k=1}^{n+m} \sigma_k^i(Y_t) \, dW_t^k \right\}, \qquad S^i(0) = s^i, i = 1, \ldots, m,$$

where $W_t = (W_t^k)_{k=1,\ldots,(n+m)}$ is an $(n+m)$-dimensional standard Brownian motion process defined on a filtered probability space $(\Omega, \mathcal{F}, P, \mathcal{F}_t)$. The factor process $Y_t$ with $n$ components is described by

$$(2.3) \quad dY_t = \beta(Y_t) \, dt + \lambda(Y_t) \, dW_t, \qquad Y(0) = y \in \mathbb{R}^n.$$

In this paper we assume that $r(y), \alpha(y), \sigma(y), \beta(y)$ and $\lambda(y)$ are given by

$$(2.4) \quad \begin{aligned} r(y) &:= r, \\ \alpha(y) &:= a + Ay, \qquad \sigma(y) := \Sigma, \\ \beta(y) &:= b + By \quad \text{and} \quad \lambda(y) := \Lambda \end{aligned}$$

with constants $r \geq 0$, $a \in \mathbb{R}^m$ and $b \in \mathbb{R}^n$. Moreover, $A, B, \Sigma, \Lambda$ are $m \times n$, $n \times n$, $m \times (n+m)$, $n \times (n+m)$ constant matrices, respectively. In this paper, we assume the following conditions (A):

(A1) $\Sigma\Sigma^* > 0,$

(A2) $G := B - \Lambda\Sigma^*(\Sigma\Sigma^*)^{-1}A$ is stable.

REMARK 2.1. The matrix $G$ is stable if the real parts of its eigenvalues are negative.

Consider an investor who invests at time $t$ a proportion $\pi_t^i$ of his wealth in the $i$th risky stock $S^i, i = 1, \ldots, m$. With $\pi_t = (\pi_t^1, \ldots, \pi_t^m)^*$ chosen, the proportion of wealth invested in the bank account is $1 - \pi_t^* \mathbf{1}$ where $(\cdot)^*$ denotes the transpose of a vector or a matrix. Here $\mathbf{1} = (1, \ldots, 1)^*$. We allow short selling ($\pi_t^i < 0$ for some $i$) or borrowing ($1 - \pi_t^* \mathbf{1} < 0$).

We assume the self-financing condition. Then the investor's wealth, $X_t^\pi$, starting with the initial capital $x$ satisfies the equation

$$(2.5) \qquad \begin{cases} \dfrac{dX_t^\pi}{X_t^\pi} = (1 - \pi_t^* \mathbf{1}) \dfrac{dS_t^0}{S_t^0} + \sum_{i=1}^m \pi_t^i \dfrac{dS_t^i}{S_t^i}, \\ X_0^\pi = x. \end{cases}$$

This equation can be solved and we have

$$(2.6) \qquad \begin{aligned} X_T^\pi = x \exp \Bigg[ \int_0^T \Big\{ r + (a + AY_t - r\mathbf{1})^* \pi_t \\ - \frac{1}{2} \pi_t^* \Sigma\Sigma^* \pi_t \Big\} dt + \int_0^T \pi_t^* \Sigma \, dW_t \Bigg]. \end{aligned}$$

For a finite $T$ and given target growth rate $k$, we shall consider the probability of minimizing the "down-side" risk,

$$(2.7) \qquad \inf_{\pi \in \mathcal{A}_T} P\left( \frac{\log X_T^\pi}{T} \le k \right),$$

where $\mathcal{A}_T$ is the set of all admissible investment strategies which will be prescribed in Section 2 [see (2.31)]. Here, the investor is interested in minimizing the probability that his wealth falls below a target growth rate. We will be mainly concerned with the large $T$ asymptotics. That is,

$$(2.8) \qquad \Pi(k) := \varlimsup_{T \to \infty} \frac{1}{T} \inf_{\pi \in \mathcal{A}_T} \log P\left( \frac{\log X_T^\pi}{T} \le k \right).$$

We shall calculate $\Pi(k)$. Connecting to this, we will obtain a (nearly) optimal strategy for (2.7) for large $T$.

2.2. *Portfolio optimization problem.* Since it is difficult to directly calculate (2.7) and (2.8), in order to solve our problem we need to introduce the following portfolio optimization problem of maximizing the expected utility. For $\gamma < 0$, we consider

$$\sup_{\pi \in \mathcal{A}_T} E\left( \frac{1}{\gamma} (X_T^\pi)^\gamma \right).$$

This is the same as

$$(2.9) \qquad \inf_{\pi \in \mathcal{A}_T} E((X_T^\pi)^\gamma).$$

The large $T$ asymptotics are given by

$$(2.10) \qquad \Psi(\gamma) := \varliminf_{T \to \infty} \frac{1}{T} \inf_{\pi \in \mathcal{A}_T} \log E(X_T^\pi)^\gamma, \qquad \gamma < 0.$$

For (2.9), we define

$$(2.11) \qquad J(x, y, \pi; T) := x^{-\gamma} E(X_T^\pi)^\gamma,$$

where $X_0^\pi = x, Y_0 = y$, and $\pi$ is taken from the set $\mathcal{A}_T$. $\mathcal{A}_T$ will be defined in (2.31). We then define

$$(2.12) \qquad v(t, y) = \inf_{\pi \in \mathcal{A}_{T-t}} \log J(x, y, \pi; T - t), \qquad \gamma < 0.$$

Using (2.6), we have

$$J(x, y, \pi; T) = E\left( \exp\left( \int_0^T \gamma \phi(Y_t, \pi_t) \, dt \right) \zeta_T^\pi \right),$$

where $\phi(y, p)$ is defined by

$$(2.13) \qquad \phi(y, p) := r + (\alpha(y) - r\mathbf{1})^* p - \frac{1 - \gamma}{2} p^* \Sigma \Sigma^* p$$

and

$$(2.14) \qquad \zeta_T^\pi = \exp\left( \int_0^T \gamma \pi_t^* \Sigma \, dW_t - \frac{1}{2} \int_0^T \gamma^2 |\Sigma^* \pi_t|^2 \, dt \right).$$

If we assume

$$(2.15) \qquad E(\zeta_T^\pi) = 1,$$

then by Girsanov theorem, we have

$$(2.16) \qquad J(x, y, \pi; T) = E^\pi \left( \exp\left( \int_0^T \gamma \phi(Y_t, \pi_t) \, dt \right) \right),$$

where $E^\pi(\cdot)$ is the expectation with respect to the probability measure $P^\pi$ defined by

$$(2.17) \qquad \left. \frac{dP^\pi}{dP} \right|_{\mathcal{F}_T} = \zeta_T^\pi.$$

We can rewrite the equation for $Y_t$ as

$$(2.18) \qquad dY_t = \beta(Y_t, \pi_t) \, dt + \Lambda \, dW^\pi(t),$$

where

$$(2.19) \qquad \beta(y, \pi) = \beta(y) + \gamma \Lambda \Sigma^* \pi$$

and

$$(2.20) \qquad W_t^\pi = W_t - \int_0^t \gamma \Sigma^* \pi_s \, ds.$$

$W_t^\pi$ is a $(n + m)$-dimensional Brownian motion under $P^\pi$. (2.9) becomes a stochastic control problem with criterion $\log J(x, y, \pi; T)$ [given in (2.16)], (2.18) is the state dynamic and $\pi_t$ is the control process. Here we note that $J(x, y, \pi; T)$ is not dependent on $x$. Therefore the value function is given by $v(0, y)$ [and $v(t, y)$ for the problem in (2.12)]. Then, by Bellman's dynamic programming principle, $v$ should satisfy the following Bellman equation:

$$\begin{cases} \dfrac{\partial v}{\partial t} + \inf_{\pi \in \mathbb{R}} \left\{ \dfrac{1}{2} \operatorname{tr}(\Lambda\Lambda^* D^2 v) + (\beta + \gamma\Lambda\Sigma^*\pi)^* Dv \right. \\ \left. \qquad\qquad + \dfrac{1}{2}(Dv)^* \Lambda\Lambda^* Dv + \gamma\phi(\cdot, \pi) \right\} = 0, \\ v(T, y) = 0, \end{cases}$$

or, equivalently,

$$(2.21) \quad \begin{cases} \dfrac{\partial v}{\partial t} + \dfrac{1}{2}\operatorname{tr}(\Lambda\Lambda^* D^2 v) + \left\{ \beta + \dfrac{\gamma}{1-\gamma}\Lambda\Sigma^*(\Sigma\Sigma^*)^{-1}(\alpha - r\mathbf{1}) \right\}^* Dv \\ \qquad + \dfrac{1}{2}(Dv)^* \Lambda \left\{ I + \dfrac{\gamma}{1-\gamma}\Sigma^*(\Sigma\Sigma^*)^{-1}\Sigma \right\} \Lambda^* Dv \\ \qquad + \dfrac{\gamma}{2(1-\gamma)}(\alpha - r\mathbf{1})^*(\Sigma\Sigma^*)^{-1}(\alpha - r\mathbf{1}) + \gamma r = 0, \\ v(T, y) = 0 \end{cases}$$

(see Section VI.8 [17]). Actually, given a solution $v(t, y)$ of (2.21), under some suitable conditions, $\widehat{\pi}_t$ defined below gives an optimal strategy of (2.9):

$$(2.22) \qquad \begin{aligned} \widehat{\pi}_t &:= \widehat{\pi}(t, Y_t), \\ \widehat{\pi}(t, y) &:= \frac{1}{1-\gamma}(\Sigma\Sigma^*)(y)^{-1}(\alpha(y) - r\mathbf{1} + \Sigma\Lambda^* Dv(t, y)) \end{aligned}$$

(see [17, 28]). Moreover, in relation to (2.10) we shall consider the ergodic type Bellman equation which is the limiting equation of the above Bellman equation,

$$\chi = \inf_{\pi \in \mathbb{R}} \left\{ \frac{1}{2}\operatorname{tr}(\Lambda\Lambda^* D^2 \xi) + (\beta + \gamma\Lambda\Sigma^*\pi)^* D\xi \right. \\ \left. \qquad\qquad + \frac{1}{2}(D\xi)^* \Lambda\Lambda^* D\xi + \gamma\phi(\cdot, \pi) \right\},$$

or, equivalently,

$$\chi = \frac{1}{2}\operatorname{tr}(\Lambda\Lambda^* D^2 \xi) + \left\{\beta + \frac{\gamma}{1-\gamma}\Lambda\Sigma^*(\Sigma\Sigma^*)^{-1}(\alpha - r\mathbf{1})\right\}^* D\xi$$

$$(2.23) \qquad + \frac{1}{2}(D\xi)^*\Lambda\left\{I + \frac{\gamma}{1-\gamma}\Sigma^*(\Sigma\Sigma^*)^{-1}\Sigma\right\}\Lambda^* D\xi$$

$$+ \frac{\gamma}{2(1-\gamma)}(\alpha - r\mathbf{1})^*(\Sigma\Sigma^*)^{-1}(\alpha - r\mathbf{1}) + \gamma r.$$

REMARK 2.2. The portfolio optimization problem in an infinite time horizon is closely connected with (2.10). On an infinite time horizon the criterion to be minimized is

$$\varliminf_{T\to\infty}\frac{1}{T}\log E(X_T^\pi)^\gamma,$$

$\pi \in \mathcal{A}$, $\mathcal{A}$ is a class of admissible strategies. The problem is to calculate the value

$$(2.24) \qquad \inf_{\pi\in\mathcal{A}}\varliminf_{T\to\infty}\frac{1}{T}\log E(X_T^\pi)^\gamma$$

and obtain an optimal strategy. The Bellman equation for (2.24) is also given by (2.23). Connecting with (2.24), we may consider the minimization of the "down-side" risk probabilities in an infinite time horizon,

$$(2.25) \qquad \inf_{\pi\in\mathcal{A}}\varliminf_{T\to\infty}\frac{1}{T}\log P\left(\frac{\log X_T^\pi}{T} \le k\right).$$

There are some more mathematical difficulties that arise from the problems in an infinite time horizon. Some further discussion for this problem will be given after we state our main theorem (Theorem 2.1).

2.3. *Main results.* The main idea of a stochastic control method for (2.9) or (2.10) is to solve the equations (2.21) or (2.23). From this, we can obtain the value in (2.9) or (2.10). We may also obtain an optimal strategy for (2.9) from the solution of (2.21). We state in this subsection the relevant results. Theorem 2.1 is our main result which also shows the connection between the problems (2.8) and (2.9) [or (2.10)]. Some propositions are given that will be used to prove Theorem 2.1.

PROPOSITION 2.1 [25]. *Assume (A1) and $\gamma < 0$. Then (2.21) has a solution given by*

$$(2.26) \qquad v(t,y) = \tfrac{1}{2}y^* P(t)y + q(t)^* y + h(t),$$

where $P(t)$ satisfies the Riccati differential equation:

$$(2.27) \quad \begin{cases} \dot{P}(t) + P(t)\Lambda N^{-1}\Lambda^* P(t) + K_1^* P(t) + P(t)K_1 - C^*C = 0, \\ P(T) = 0, \end{cases}$$

where

$$(2.28) \quad \begin{aligned} N^{-1} &:= \left\{ I + \frac{\gamma}{1-\gamma}\Sigma^*(\Sigma\Sigma^*)^{-1}\Sigma \right\} > 0, \\ K_1 &:= B + \frac{\gamma}{1-\gamma}\Lambda\Sigma^*(\Sigma\Sigma^*)^{-1}A, \\ C &:= \sqrt{-\frac{\gamma}{1-\gamma}}\Sigma^*(\Sigma\Sigma^*)^{-1}A, \end{aligned}$$

and $q(t), h(t)$ satisfy the equations

$$(2.29) \quad \begin{cases} \dot{q}(t) + (K_1 + \Lambda N^{-1}\Lambda^* P(t))^* q(t) + P(t)b \\ \qquad + \frac{\gamma}{1-\gamma}(A^* + P(t)\Lambda\Sigma^*)(\Sigma\Sigma^*)^{-1}(a - r\mathbf{1}) = 0, \\ q(T) = 0, \end{cases}$$

$$(2.30) \quad \begin{cases} \dot{h}(t) + \frac{1}{2}\operatorname{tr}(\Lambda\Lambda^* P(t)) + \frac{1}{2}q(t)^*\Lambda\Lambda^* q(t) + b^* q(t) + \gamma r \\ \qquad + \frac{\gamma}{2(1-\gamma)}(a - r\mathbf{1} + \Sigma\Lambda^* q(t))^*(\Sigma\Sigma^*)^{-1} \\ \qquad \times (a - r\mathbf{1} + \Sigma\Lambda^* q(t)) = 0, \\ h(T) = 0. \end{cases}$$

The proof can be found in [25]. We will give some remarks in Section 4.1. We now define the class of admissible investment strategies, $\mathcal{A}_T$,

$$(2.31) \quad \begin{aligned} \mathcal{A}_T := \Bigg\{ &(\pi_t)_{t\in[0,T]}; \\ &E\bigg[\mathcal{E}\bigg(\int [(P(s)Y_s + q(s))^*\Lambda + \gamma\pi_s^*\Sigma]\,dW_s\bigg)_T\bigg] = 1 \Bigg\}. \end{aligned}$$

This class of admissible strategies is also used in [25]. Here we use the notation, $\mathcal{E}(Z) := (\mathcal{E}(Z)_t)_{t\in[0,T]}$, for the stochastic exponential of a continuous semimartingale $Z : \mathcal{E}(Z)_t := e^{Z_t - 1/2\langle Z\rangle_t}$. Therefore, $\zeta_T^\pi$ in (2.14) is equal to $\mathcal{E}(\int\gamma\pi^*\Sigma\,dW)_T$.

PROPOSITION 2.2 [25].   *Assume (A1) and $\gamma < 0$. Let $P(t), q(t), h(t)$ be defined as in Proposition 2.1. Then the following defines an admissible strategy:*

$$\hat{\pi}(t) := \hat{\pi}(t, Y_t),$$

(2.32)
$$\widehat{\pi}(t, y) = \frac{1}{1 - \gamma}(\Sigma\Sigma^*)^{-1}[a - r\mathbf{1} + \Sigma\Lambda^* q(t) + \{A + \Sigma\Lambda^* P(t)\}y]$$

and is optimal for the problem (2.9). Moreover,

$$v(0, y) = \tfrac{1}{2}y^* P(0)y + q(0)^* y + h(0).$$

The proofs of Proposition 2.2 can be found in [25]. Some remarks will also be given in Section 4.1.

In the following, notation $v(t, y; T), P(t; T), q(t; T), h(t; T)$ will be used for the dependence of $P(t), q(t), h(t)$ on $T$. Similarly, we use

$$v(t, y; T; \gamma), \qquad P(t; T; \gamma), \qquad q(t; T; \gamma), \qquad h(t; T; \gamma),$$

if we need to discuss the dependence of the functions on $\gamma$.

We now consider (2.10). We consider a solution of (2.23) that is quadratic in $y$. That is,

(2.33)
$$\xi(y) := \tfrac{1}{2}y^* \overline{P}y + \overline{q}^* y,$$

where $\overline{P}$ is a symmetric $n \times n$ matrix, and $\overline{q}$ is a vector in $\mathbb{R}^n$. $\overline{P}$ and $\overline{q}$ will satisfy the equations:

(2.34)
$$K_1^* \overline{P} + \overline{P}K_1 + \overline{P}\Lambda N^{-1}\Lambda^* \overline{P} - C^*C = 0,$$

(2.35)
$$(K_1 + \Lambda N^{-1}\Lambda^* \overline{P})^* \overline{q} + \overline{P}b$$
$$+ \frac{\gamma}{1 - \gamma}(A^* + \overline{P}\Lambda\Sigma^*)(\Sigma\Sigma^*)^{-1}(a - r\mathbf{1}) = 0.$$

We have the following results.

PROPOSITION 2.3 [25]. *In addition to* (A1), *we assume* (A2). *Then the following properties hold:*

(i) $P(t) = P(t; T)$ *converges as* $T \to \infty$ *to a nonpositive definite matrix* $\overline{P}$, *which is a solution of the algebraic Riccati equation* (2.34). *Moreover,*

$$K_1 + \Lambda N^{-1}\Lambda^* \overline{P}$$

*is stable, and* $\overline{P}$ *satisfies the estimate*

(2.36)
$$-\int_0^\infty e^{sG^*}A^*(\Sigma\Sigma^*)^{-1}Ae^{sG}\,ds \leq \overline{P} \leq 0,$$

*where* $G$ *is given in* (A2).

(ii) $q(t) = q(t; T)$ *converges as* $T \to \infty$ *to a constant vector* $\overline{q}$ *which satisfies (2.35). Moreover,* $-\dot{h}(t) = -\dot{h}(t; T)$ *converges to a constant* $\chi(\gamma)$ *defined by*

$$
(2.37) \quad \begin{aligned}
\chi(\gamma) = &\frac{1}{2} \operatorname{tr}(\Lambda \Lambda^* \overline{P}) + \frac{1}{2} \overline{q}^* \Lambda \Lambda^* \overline{q} + b^* \overline{q} + \gamma r \\
&+ \frac{\gamma}{2(1-\gamma)} (a - r\mathbf{1} + \Sigma \Lambda^* \overline{q})^* (\Sigma \Sigma^*)^{-1} (a - r\mathbf{1} + \Sigma \Lambda^* \overline{q}).
\end{aligned}
$$

(iii) *We obtain*

$$
(2.38) \quad \lim_{T \to \infty} \frac{v(0; y)}{T} = \lim_{T \to \infty} \frac{h(0; T)}{T} = \chi(\gamma).
$$

The proof of Proposition 2.3 can be found in [25]. Some remarks will be also given in Section 4.2.

We use $\xi_\gamma(y), \overline{P}(\gamma), \overline{q}(\gamma)$ for the dependence of $\xi(y), \overline{P}, \overline{q}$ on $\gamma$.

PROPOSITION 2.4.    *Along with (A1) and (A2), we assume (A3);*

(A3)                          $(B, \Lambda)$ *is controllable.*

*Then the pair* $(\xi, \chi(\gamma))$ *is the unique solution of (2.23) with* $\chi = \chi(\gamma)$ *where* $\xi(y)$ *and* $\chi(\gamma)$ *are given by (2.33) and (2.37), respectively.*

The proof of Proposition 2.4 will be given in Section 4.2.

REMARK 2.3.    The pair $(K, L)$, of the $n \times n$ matrix $K$ and the $n \times l$ matrix $L$, is said controllable if the $n \times nl$ matrix $(L, KL, K^2L, \ldots, K^{n-1}L)$ has rank $n$. We remark that the generator of $Y_t$, $\mathcal{G}f = \frac{1}{2} \operatorname{tr}(\Lambda \Lambda^* D^2 f) + \beta \cdot Df$, is an hypoelliptic second-order operator if (A3) holds (see [20]).

REMARK 2.4.    It is shown in [22] that there is $\chi^*(\gamma)$ such that for any $\chi \geq \chi^*(\gamma)$, (2.23) has a solution $\xi(y)$. There are only finitely many particular $\chi$'s where the solution $\xi(y)$ is quadratic in $y$. Under certain conditions, (2.23), for $\chi = \chi^*(\gamma)$, has a unique solution $\xi^*(y)$ satisfying $\xi^*(0) = 0$. The unique pair $\chi^*(\gamma), \xi^*(\cdot)$ may also be characterized by the growth condition of $\xi^*(y)$ as $|y| \to \infty$ (see [27]).

The differentiability of $\xi_\gamma(y)$ [or $\chi(\gamma), \overline{P}(\gamma), \overline{q}(\gamma)$] will play an important role in the proof of our main theorem (Theorem 2.1). We have the following results.

PROPOSITION 2.5.    *Under assumptions (A1) and (A2), the following results hold:*

(i) $\chi(\gamma)$ and $\overline{P}(\gamma), \overline{q}(\gamma)$ are twice differentiable with respect to $\gamma$.

(ii) $\chi(\gamma)$ is convex with respect to $\gamma$.

The proof of Proposition 2.5(i) [and (ii)] will be given in Section 4.3 (and Section 4.4).

Similar to (2.32), we use $\xi(\cdot)$ to define

$$
\begin{aligned}
(2.39) \qquad \pi_\gamma(y) &= \frac{1}{1-\gamma}(\Sigma\Sigma^*)^{-1}(\alpha(y) - r\mathbf{1} + \Sigma\Lambda^* D\xi(y)) \\
&= \frac{1}{1-\gamma}(\Sigma\Sigma^*)^{-1}[a - r\mathbf{1} + \Sigma\Lambda^*\overline{q}(\gamma) + \{A + \Sigma\Lambda^*\overline{P}(\gamma)\}y].
\end{aligned}
$$

Now we also define

$$
\begin{aligned}
(2.40) \qquad \beta_\gamma(y) &:= \beta(y) + \gamma\Lambda\Sigma^*\pi_\gamma(y) + \Lambda\Lambda^* D\xi(y) \\
&= b + By + \frac{\gamma}{1-\gamma}\Lambda\Sigma^*(\Sigma\Sigma^*)^{-1}(a + Ay - r\mathbf{1}) \\
&\quad + \Lambda N^{-1}\Lambda^*(\overline{P}(\gamma)y + \overline{q}(\gamma)) \\
&= (K_1 + \Lambda N^{-1}\Lambda^*\overline{P}(\gamma))y + f_\gamma,
\end{aligned}
$$

where

$$
f_\gamma := b + \frac{\gamma}{1-\gamma}\Lambda\Sigma^*(\Sigma\Sigma^*)^{-1}(a - r\mathbf{1}) + \Lambda N^{-1}\Lambda^*\overline{q}.
$$

Define

$$
\begin{aligned}
(2.41) \qquad u = u(y) &:= \{\Lambda^* D\xi(y) + \gamma\Sigma^*\pi_\gamma(y)\}^* \\
&= \Bigg\{ \Lambda^*(\overline{P}y + \overline{q}) \\
&\quad + \frac{\gamma}{1-\gamma}\Sigma^*(\Sigma\Sigma^*)^{-1}(\Sigma\Lambda^*(\overline{P}y + \overline{q}) + a + Ay - r\mathbf{1}) \Bigg\}^*.
\end{aligned}
$$

We will show

$$
(2.42) \qquad E\left(\mathcal{E}\left(\int u\, dW\right)_t\right) = 1, \qquad t > 0.
$$

Then we define a new probability measure $\widehat{P}$ by

$$
(2.43) \qquad \frac{d\widehat{P}}{dP}\Bigg|_{\mathcal{F}_t} := \mathcal{E}\left(\int u\, dW\right)_t.
$$

Define $\widehat{W}_t$ by

$$
(2.44) \qquad \widehat{W}_t = W_t - \int_0^t u(Y_s)^*\, ds.
$$

Then $\widehat{W}_t$ is a Brownian motion under the probability measure $\widehat{P}$, and $Y$ satisfies

$$(2.45) \qquad\qquad dY_t = \beta_\gamma(Y_t)\,dt + \Lambda\,d\widehat{W}_t.$$

REMARK 2.5.   We compare $\beta_\gamma(y)$ in (2.40) and $\beta(\cdot, \cdot)$ in (2.19),

$$\beta_\gamma(y) = \beta(y, \pi_\gamma(y)) + \Lambda\Lambda^*(\overline{P}(\gamma)y + \overline{q}(\gamma)).$$

As shown in (2.16), $\beta(y, \pi_\gamma(y))$ is used to change the measure to derive the useful expression of $J(x, y, \pi; T)$. Because of the integral from 0 to $T$ in (2.16) when $\pi_t = \pi_\gamma(Y_t)$, the expectation may grow exponentially and will cause difficulty in the analysis. The difference of $\beta_\gamma(y)$ and $\beta(y, \pi_\gamma(y))$ is made to take care of this. This gives another measure change, hence another term is added to $\beta(y, \pi_\gamma(y))$ which leads to $\beta_\gamma(y)$. A similar idea is used in [23]. In the proof of our main theorem (Theorem 2.1), we will also see some uses of $\widehat{P}$.

Under $\widehat{P}$, $Y_t$ is a Gaussian process. The variance of $Y_t$ is given by $U(t)$,

$$U(t) = \int_0^t e^{(t-s)(K_1 + \Lambda N^{-1}\Lambda^*\overline{P})}\Lambda\Lambda^* e^{(t-s)(K_1 + \Lambda N^{-1}\Lambda^*\overline{P})^*}\,ds,$$

and its mean $m(t)$ is the solution of the following equation:

$$\dot{m}(t) = (K_1 + \Lambda N^{-1}\Lambda^*\overline{P})m(t) + f_\gamma.$$

We show in the next proposition that $K_1 + \Lambda^* N^{-1}\Lambda^*\overline{P}$ is stable under assumption (A2); then $Y_t$ is ergodic under $\widehat{P}$.

PROPOSITION 2.6.   *Under assumptions (A1) and (A2), $K_1 + \Lambda^* N^{-1}\Lambda^*\overline{P}$ is stable. There are $c_1(\gamma) > 0, c_2(\gamma) > 0$ and a positive definite $n \times n$ matrix $K_\gamma$ such that*

$$(2.46) \qquad\qquad y^* K_\gamma \beta_\gamma(y) \le -c_1(\gamma)|y|^2 + c_2(\gamma),$$

*where $\beta_\gamma(y)$ is given by (2.40).*

The proof is given in Section 4.5.

PROPOSITION 2.7.   *Under assumptions (A1) and (A2) we obtain*

$$(2.47) \qquad\qquad \lim_{\gamma \to -\infty} \chi'(\gamma) = r.$$

For the proof, see Section 4.6. This result gives the range of $k$ in (2.8) that our main result holds (Theorem 2.1).

We can now state our main theorem. Its proof, given in next section, is based on the results in Propositions 2.1–2.7.

THEOREM 2.1. *Assume (A1) and (A2). Then, for all $r < k < \chi'(0-)$, we have the following:*

$$\Pi(k) = \inf_{\gamma < 0} \{\chi(\gamma) - \gamma k\}. \tag{2.48}$$

*Moreover, define the strategy,*

$$\widehat{\pi}_t^{[k]} := \widehat{\pi}(t, Y_t), \tag{2.49}$$

*where $\widehat{\pi}(t, y)$ is defined in (2.32) with $\gamma = \gamma(k)$. Here $\gamma(k) < 0$ s.t. $\chi'(\gamma(k)) = k \in (r, \chi'(0-))$. We denote $\chi'(\gamma) = \frac{d\chi}{d\gamma}(\gamma)$. Then $\widehat{\pi}^{[k]}$, on the sufficiently large time horizon $T$, is a nearly optimal strategy for the problem (2.7), namely,*

$$\Pi(k) = \lim_{T \to \infty} \frac{1}{T} \log P\left(\frac{\log X_T^{\widehat{\pi}^{[k]}}}{T} \le k\right).$$

*For $k < r$,*

$$\Pi(k) = \inf_{\gamma < 0} \{\chi(\gamma) - \gamma k\} = -\infty.$$

*If $B$ is stable, then*

$$\chi(0-) = 0 \tag{2.50}$$

*and*

$$\chi'(0-) = \frac{1}{2} \operatorname{tr}\left(\Lambda\Lambda^* \frac{d\overline{P}}{d\gamma}(0-)\right) \tag{2.51}$$

$$+ \frac{1}{2}(AB^{-1}b - (a - r\mathbf{1}))^* (\Sigma\Sigma^*)^{-1}(AB^{-1}b - (a - r\mathbf{1})) + r,$$

*where*

$$\frac{d\overline{P}}{d\gamma}(0-) = \int_0^\infty e^{sB^*} A^*(\Sigma\Sigma^*)^{-1} A e^{sB} \, ds,$$

$$\Pi(k) = \inf_{\gamma < 0} \{\chi(\gamma) - \gamma k\} = 0, \qquad k \ge \chi'(0-).$$

REMARK 2.6. In Proposition 2.2, for a given $T$,

$$P(t) = P(t; T), \qquad q(t) = q(t; T), \qquad h(t) = h(t; T)$$

and

$$\widehat{\pi}(t, y) = \widehat{\pi}(t, y; T).$$

Hence, to avoid confusion, we may use $\widehat{\pi}_t^{T,k}$ for $\pi_t^{[k]}$. Therefore, for $\gamma = \gamma(k)$,

$$E((X_T^{\widehat{\pi}^{T,k}})^\gamma) = \inf_{\pi \in \mathcal{A}_T} E((X_T^\pi)^\gamma).$$

$\widehat{\pi}_t^{T,k}$ is an optimal strategy for (2.9) $[\gamma = \gamma(k)]$ but may not be an optimal strategy for (2.7). This is the reason we say $\widehat{\pi}_t^{T,k}$ is nearly optimal for (2.7) when $T$ is large, since the value using this strategy is close to the optimal value when $T$ is large.

On the other hand, we can use $\pi_\gamma(\cdot)$ in (2.39) to define

$$\widehat{\pi}_\gamma(t) = \pi_\gamma(Y_t).$$

We may expect that this will also give a nearly optimal strategy. That is, in a sense,

(2.52)                          $$E((X_T^{\widehat{\pi}_\gamma})^\gamma)$$

is close to

$$\inf_{\pi \in \mathcal{A}_T} E((X_T^\pi)^\gamma),$$

if $T$ is large. There are two problems when one wants to prove this rigorously. For the first problem, it is not easy to show $\widehat{\pi}_\gamma(\cdot)$ is an element of $\mathcal{A}_T$. For the second problem, it may happen that (2.52) becomes infinite for some finite $T$. See [14] for a study of a model where there is one stock in the market. When this happens, the problem (2.24) may not have a solution. That is, there is no optimal strategy for (2.24). However, it is shown in [14] that some modification of $\widehat{\pi}_\gamma(\cdot)$ gives a nearly optimal strategy. Such behavior also indicates that (2.25) may be more difficult to treat than (2.8). However, if we assume (A1)–(A3) and the following condition,

$$\overline{P}(\gamma)\Lambda\Sigma^*(\Sigma\Sigma^*)^{-1}\Sigma\Lambda^*\overline{P}(\gamma) < A^*(\Sigma\Sigma^*)^{-1}A,$$

then it is proved in [25] (Theorem 2.3) that

$$\lim_{T \to \infty} \frac{1}{T} \log E(X_T^{\widehat{\pi}_\gamma}) = \chi(\gamma).$$

Using this and following the same argument as in Theorem 2.1, we can obtain the upper estimate for down-side risk probability. For the lower estimate, the same argument in Theorem 2.1 can be applied. In conclusion, we have the following result; its proof is omitted.

THEOREM 2.2.  *Assume (A1)–(A3). Let $r < k < \chi'(0-)$ and $\gamma(k) < 0$ be the unique number satisfying $\chi'(\gamma(k)) = k$. Assume*

$$\overline{P}(\gamma(k))\Lambda\Sigma^*(\Sigma\Sigma^*)^{-1}\Sigma\Lambda^*\overline{P}(\gamma(k)) < A^*(\Sigma\Sigma^*)^{-1}A.$$

*Define*

$$\widehat{\pi}_k(t) := \pi_{\gamma(k)}(Y_t).$$

*Then $\widehat{\pi}_k(\cdot) \in \mathcal{A}_T$ for any $T$ and $\widehat{\pi}_k$ on a sufficiently large time horizon $T$ is a nearly optimal strategy for the problem (2.7). Namely,*

$$\Pi(k) = \lim_{T\to\infty} \frac{1}{T} \log P\left(\frac{\log X_T^{\widehat{\pi}_k}}{T} \le k\right).$$

*In (2.25) we define $\mathcal{A}$ to be a family consisting of processes $\pi_t$ such that $\pi|_{[0,T]}$, the restriction of $\pi_t$ to $[0,T]$, is in $\mathcal{A}_T$ for all $T$. Then $\widehat{\pi}_k(\cdot) \in \mathcal{A}$ and is an optimal strategy for (2.25).*

**3. Proof of Theorem 2.1.** In this section, we shall give a proof of our main theorem using results in Propositions 2.1–2.7. The proof of the propositions will be given in Section 4.

From Proposition 2.5(ii), $\chi$ is convex. Let us consider

$$\widehat{\chi}(k) := \inf_{\gamma<0}\{\chi(\gamma) - \gamma k\}, \qquad k > r.$$

Since $\chi$ is smooth, we see that $\widehat{\chi}(k)$ is strictly concave, nondecreasing and satisfies

$$\widehat{\chi}(k) = \begin{cases} 0, & \text{if } k \ge \chi'(0-), \\ \chi(\gamma(k)) - \gamma(k)\chi'(\gamma(k)), & \text{if } r < k < \chi'(0-), \end{cases}$$

where $\gamma(d) < 0$ is such that $\chi'(\gamma(d)) = d \in (r, \chi'(0-))$. In Section 4.4, we will show $\chi'(\gamma) > 0$ for $\gamma < 0$. Therefore, $\gamma(d)$ is uniquely defined. Moreover, $\widehat{\chi}$ is continuous on $(r, \chi'(0-))$.

We now take a small $\varepsilon > 0$. Then $k - \varepsilon > r$. Suppose $\tilde{\gamma} < 0$ attains the following:

$$\widehat{\chi}(k - \varepsilon) := \inf_{\gamma<0}\{\chi(\gamma) - \gamma(k-\varepsilon)\} = \chi(\tilde{\gamma}) - \tilde{\gamma}(k-\varepsilon).$$

Note that

(3.1) $$k - \varepsilon = \chi'(\tilde{\gamma}).$$

We denote $\tilde{\gamma}$ as $\gamma$ in the following.

From Proposition 2.3, we recall that the pair $(\chi(\gamma), \xi)$ solves (2.23);

$$\chi(\gamma) = \frac{1}{2}\operatorname{tr}(\Lambda\Lambda^* D^2\xi) + \left\{b + By + \frac{\gamma}{1-\gamma}\Lambda\Sigma^*(\Sigma\Sigma^*)^{-1}(a + Ay - r\mathbf{1})\right\}^* D\xi$$

$$+ \frac{1}{2}(D\xi)^*\Lambda N^{-1}\Lambda^* D\xi$$

$$+ \frac{\gamma}{2(1-\gamma)}(a + Ay - r\mathbf{1})^*(\Sigma\Sigma^*)^{-1}(a + Ay - r\mathbf{1}) + \gamma r,$$

or, equivalently,

$$(3.2) \qquad \chi(\gamma) = \tfrac{1}{2}\operatorname{tr}(\Lambda\Lambda^* D^2\xi(y)) + \beta_\gamma(y)^* D\xi(y) + V_0,$$

where $\beta_\gamma(y)$ is defined in (2.40) and

$$
\begin{aligned}
V_0 &= V_0(y) \\
(3.3) \qquad &:= -\tfrac{1}{2}D\xi(y)^*\Lambda N^{-1}\Lambda D\xi(y) \\
&\quad + \frac{\gamma}{2(1-\gamma)}(a + Ay - r\mathbf{1})^*(\Sigma\Sigma^*)^{-1}(a + Ay - r\mathbf{1}) + \gamma r.
\end{aligned}
$$

By differentiation of (3.2) with respect to $\gamma$, we have

$$
\begin{aligned}
\chi'(\gamma) &= \tfrac{1}{2}\operatorname{tr}(\Lambda\Lambda^* D^2\eta(y)) + \beta_\gamma(y)^* D\eta(y) \\
&\quad + \frac{1}{(1-\gamma)^2}\Lambda\Sigma^*(\Sigma\Sigma^*)^{-1}(a + Ay - r\mathbf{1})^* D\xi \\
&\quad + \frac{1}{2(1-\gamma)^2}(D\xi)^*(y)\Lambda\Sigma^*(\Sigma\Sigma^*)^{-1}\Sigma\Lambda^* D\xi(y) \\
&\quad + \frac{1}{2(1-\gamma)^2}(a + Ay - r\mathbf{1})^*(\Sigma\Sigma^*)^{-1}(a + Ay - r\mathbf{1}) + r,
\end{aligned}
$$

where $\eta = \frac{\partial\xi}{\partial\gamma}$. That is,

$$(3.4) \quad \chi'(\gamma) = \tfrac{1}{2}\operatorname{tr}(\Lambda\Lambda^* D^2\eta(y)) + \beta_\gamma(y)^* D\eta(y) + V_1,$$

$$
\begin{aligned}
V_1 = V_1(y) &:= \frac{1}{2(1-\gamma)^2}\{\Sigma\Lambda^* D\xi(y) + (a + Ay - r\mathbf{1})\}^*(\Sigma\Sigma^*)^{-1} \\
(3.5) \qquad &\quad \times \{\Sigma\Lambda^* D\xi(y) + (a + Ay - r\mathbf{1})\} + r.
\end{aligned}
$$

Furthermore, from (3.2) and (3.4), we can obtain

$$
\begin{aligned}
(3.6) \qquad &\chi(\gamma) - \gamma\chi'(\gamma) \\
&= \tfrac{1}{2}\operatorname{tr}(\Lambda\Lambda^* D^2(\xi - \gamma\eta)(y)) + \beta_\gamma(y)^* D(\xi - \gamma\eta)(y) + V_2,
\end{aligned}
$$

where $V_2$ is defined by

$$(3.7) \qquad\qquad V_2 = V_2(y) = -\tfrac{1}{2}|u|^2;$$

$u(\cdot)$ is in (2.41). Indeed, from (3.2) and (3.4),

$$\chi(\gamma) - \gamma\chi'(\gamma) = \tfrac{1}{2}\operatorname{tr}(\Lambda\Lambda^* D^2(\xi - \gamma\eta)(y)) + \beta_\gamma(y)^* D(\xi - \gamma\eta)(y) + V_0 - \gamma V_1.$$

By (3.3) and (3.5), a straightforward calculation shows that

$$
\begin{aligned}
V_0 - \gamma V_1 = -\frac{1}{2} & D\xi(y)^* \Lambda \Lambda^* D\xi(y) \\
& - \frac{\gamma}{1-\gamma} \{\Sigma \Lambda^* D\xi(y) + (a + Ay - r\mathbf{1})\}^* (\Sigma\Sigma^*)^{-1} \Sigma \Lambda^* D\xi(y) \\
& - \frac{1}{2} \left(\frac{\gamma}{1-\gamma}\right)^2 \{\Sigma \Lambda^* D\xi(y) + (a + Ay - r\mathbf{1})\}^* \\
& \times (\Sigma\Sigma^*)^{-1} \{\Sigma \Lambda^* D\xi(y) + (a + Ay - r\mathbf{1})\} \\
= -\frac{1}{2} & |u|^2,
\end{aligned}
$$

where the last equality follows from (2.41). Hence we have (3.7).

For $\pi \in \mathcal{A}_T$, we have by (2.6),

$$
X_T^\pi = x \exp\left[\int_0^T \left\{r + (a + AY_t - r\mathbf{1})^* \pi_t - \frac{1}{2}\pi_t^* \Sigma\Sigma^* \pi_t\right\} dt + \int_0^T \pi_t^* \Sigma \, dW_t\right].
$$

Then

$$
\begin{aligned}
\frac{\log X_T^\pi}{T} = \frac{\log x}{T} & + \frac{1}{T}\int_0^T \pi_t^* \Sigma \, dW_t \\
& + \frac{1}{T}\int_0^T \left[-\frac{1}{2}\pi_t^* \Sigma\Sigma^* \pi_t + \{\pi_t^*(a + AY_t - r\mathbf{1}) + r\}\right] dt.
\end{aligned}
$$

Recalling $\widehat{W}_t$ defined by (2.44), we can rewrite the last relation as

$$
\begin{aligned}
\frac{\log X_T^\pi}{T} = \frac{\log x}{T} & + \frac{1}{T}\int_0^T \pi_t^* \Sigma \, d\widehat{W}_t \\
& - \frac{1}{2T}\int_0^T \left\{\pi_t - \frac{1}{1-\gamma}(\Sigma\Sigma^*)^{-1}(\Sigma \Lambda^* D\xi(Y_t) + a + AY_t - r\mathbf{1})\right\}^* \\
& \times \Sigma\Sigma^* \left\{\pi_t - \frac{1}{1-\gamma}(\Sigma\Sigma^*)^{-1}(\Sigma \Lambda^* D\xi(Y_t) + a + AY_t - r\mathbf{1})\right\} dt \\
& + \frac{1}{T}\int_0^T \left\{\frac{1}{2(1-\gamma)^2}(\Sigma \Lambda^* D\xi(Y_t) + a + AY_t - r\mathbf{1})^* \right. \\
& \left. \times (\Sigma\Sigma^*)^{-1}(\Sigma \Lambda^* D\xi(Y_t) + a + AY_t - r\mathbf{1}) + r\right\} dt.
\end{aligned}
$$

From (3.5) and with some calculation, we can rewrite this as

$$
\frac{\log X_T^\pi}{T} = \frac{\log x}{T} + \frac{1}{T}\int_0^T \frac{1}{1-\gamma}\{(\Sigma\Sigma^*)^{-1}(\Sigma \Lambda^* D\xi(Y_t)
$$

$$+ a + AY_t - r\mathbf{1})\}^* \Sigma \, d\widehat{W}_t$$

(3.8)
$$+ \frac{1}{T} \log \mathcal{E}\bigg(\int \Big\{\pi - \frac{1}{1-\gamma}(\Sigma\Sigma^*)^{-1}(\Sigma\Lambda^* D\xi(Y)$$

$$+ a + AY - r\mathbf{1})\Big\}^* \Sigma \, d\widehat{W}\bigg)_T$$

$$+ \frac{1}{T} \int_0^T V_1(Y_t) \, dt.$$

Define the following events:

$$A := \left\{\omega; \frac{\log X_T^\pi}{T} \leq k\right\},$$

$$A_{1,T} := \left\{\omega; \frac{1}{T}\int_0^T V_2(Y_t) \, dt \geq \chi(\gamma) - \gamma\chi'(\gamma) - \varepsilon\right\},$$

$$A_{2,T} := \left\{\omega; \frac{1}{T}\int_0^T u(Y_t)^* \, d\widehat{W}_t \leq \varepsilon\right\}$$

and

$$A_{3,T} := \left\{\omega; \frac{\log x}{T} + \frac{1}{T}\int_0^T V_1(Y_t) \, dt \leq \chi'(\gamma) + \frac{\varepsilon}{2}\right\},$$

$$A_{4,T} := \bigg\{\omega; \frac{1}{T}\int_0^T \frac{1}{1-\gamma}\{(\Sigma\Sigma^*)^{-1}(\Sigma\Lambda^* D\xi(Y_t)$$

$$+ a + AY_t - r\mathbf{1})\}^* \Sigma \, d\widehat{W}_t \leq \frac{\varepsilon}{4}\bigg\},$$

$$A_{5,T} := \bigg\{\omega; \frac{1}{T}\log \mathcal{E}\bigg(\int\Big(\pi - \frac{1}{1-\gamma}(\Sigma\Sigma^*)^{-1}$$

$$\times (\Sigma\Lambda^* D\xi(Y) + a + AY - r\mathbf{1})\Big)^* \Sigma \, dW\bigg)_T \leq \frac{\varepsilon}{4}\bigg\}.$$

From (3.1) and (3.8), we see that

(3.9)
$$A_{3,T} \cap A_{4,T} \cap A_{5,T} \subseteq A.$$

Recall that $\widehat{P}$ is the probability defined by (2.43). By using (3.6), (2.45) and Itô's formula, we have

$$(\xi - \gamma\eta)(Y_T) - (\xi - \gamma\eta)(y) = (\chi(\gamma) - \gamma\chi'(\gamma))T - \int_0^T V_2(Y_t) \, dt$$

$$+ \int_0^T D(\xi - \gamma\eta)^*(Y_t)\Lambda \, d\widehat{W}_t.$$

We apply Chebyshev's inequality:

$$\widehat{P}(A_{1,T}^c) \le \frac{1}{\varepsilon^2} \widehat{E}\left[\left|\chi(\gamma) - \gamma\chi'(\gamma) - \frac{1}{T}\int_0^T V_2(Y_t)\,dt\right|^2\right]$$

$$\le \frac{1}{\varepsilon^2 T^2}\widehat{E}\left[\left|(\xi - \gamma\eta)(Y_T) - (\xi - \gamma\eta)(y)\right.\right.$$

$$\left.\left. - \int_0^T D(\xi - \gamma\eta)^*(Y_t)\Lambda\,d\widehat{W}_t\right|^2\right].$$

Noting that $\xi(y)$, and $\eta(y)$ are quadratic in $y$, then

$$|(\xi - \gamma\eta)(y)| \le K(1 + |y|^2),$$

$$|D(\xi - \gamma\eta)(y)| \le K(1 + |y|).$$

Therefore, we have

$$\widehat{E}\left[\left|(\xi - \gamma\eta)(Y_T) - (\xi - \gamma\eta)(y) - \int_0^T D(\xi - \gamma\eta)^*(Y_t)\lambda(Y_t)\,d\widehat{W}_t\right|^2\right]$$

$$\le K\left\{1 + |y|^4 + T + \widehat{E}[|Y_T|^4] + \int_0^T \widehat{E}[|Y_t|^2]\,dt\right\}$$

$$\le d_1(\gamma)\{|y|^4 + 1 + (|y|^2 + 1)T\}$$

for some $d_1(\gamma) > 0$ where the last inequality follows from Lemma 3.1 below. Therefore, we see that

$$\widehat{P}(A_{1,T}^c) \le \frac{d_1(\gamma)\{|y|^4 + 1 + (|y|^2 + 1)T\}}{\varepsilon^2 T^2}$$

and

$$(3.10) \qquad\qquad \widehat{P}(A_{1,T}^c) \le \varepsilon,$$

provided $T$ is sufficiently large.

By using Chebyshev's inequality again, we have

$$\widehat{P}(A_{2,T}^c) \le \frac{1}{\varepsilon^2 T^2}\widehat{E}\left[\left|\int_0^T u(Y_t)^*\,d\widehat{W}_t\right|^2\right] = \frac{1}{\varepsilon^2 T^2}\widehat{E}\left[\int_0^T |u(Y_t)|^2\,dt\right].$$

By (2.41) and using $|D\xi(y)| \le K(1 + |y|)$, we obtain

$$\widehat{P}(A_{2,T}^c) \le \frac{K}{\varepsilon^2 T^2}\left(T + \int_0^T \widehat{E}[|Y_t|^2]\,dt\right)$$

$$\le \frac{d_2(\gamma)(|y|^2 + 1)}{\varepsilon^2 T},$$

for some $d_2(\gamma) > 0$ where the last inequality follows from Lemma 3.1 below. Therefore, we see that

$$\widehat{P}(A_{2,T}^c) \le \varepsilon, \tag{3.11}$$

provided $T$ is sufficiently large.

Using a similar argument, we can show

$$\widehat{P}(A_{4,T}^c) \le \varepsilon. \tag{3.12}$$

We now consider $\widehat{P}(A_{3,T}^c)$. By (3.4) and Itô's formula,

$$d\eta(Y_t) = -V_1(Y_t)\,dt + \chi'(\gamma)\,dt + (D\eta(Y_t))^*\Lambda\,d\widehat{W}_t.$$

From this, we have

$$\int_0^T V_1(Y_t)\,dt = \eta(Y_0) - \eta(Y_T) + \chi'(\gamma)T + \int_0^T (D\eta(Y_t))^*\Lambda\,d\widehat{W}_t.$$

Then, $A_{3,T}$ is the same as

$$\left\{ \omega; \frac{\log x}{T} + \frac{\eta(Y_0) - \eta(Y_T)}{T} + \frac{1}{T}\int_0^T (D\eta(Y_t))^*\Lambda\,d\widehat{W}_t \le \frac{\varepsilon}{2} \right\}.$$

Now we can use the same argument as above to get

$$\widehat{P}(A_{3,T}^c) \le \varepsilon, \tag{3.13}$$

provided $T$ is sufficiently large.

For $A_{5,T}$, we have

$$\widehat{P}(A_{5,T}^c) \le e^{-\varepsilon/4T}$$
$$\times \widehat{E}\left[ \mathcal{E}\left( \int \left\{ \pi - \frac{1}{1-\gamma}(\Sigma\Sigma^*)^{-1}(\Sigma\Lambda^* D\xi(Y) \right.\right.\right.$$
$$\left.\left.\left. + a + AY - r\mathbf{1}) \right\}^* \Sigma\,dW \right)_T \right]$$
$$\le e^{-\varepsilon/4T}.$$

Then we have

$$\widehat{P}(A_{5,T}^c) \le \varepsilon, \tag{3.14}$$

if $T$ is sufficiently large.

Hence, from (3.7), (3.9), (3.10), (3.11), (3.12), (3.13) and (3.14), we have

$$P(A) = \widehat{E}\left[ \exp\left\{ -\int_0^T u(Y_t)^*\,d\widehat{W}_t - \frac{1}{2}\int_0^T |u(Y_t)|^2\,dt \right\}; A \right]$$

$$= \widehat{E}\Big[\exp\Big\{-\int_0^T u(Y_t)^* \, d\widehat{W}_t + \int_0^T V_2(Y_t) \, dt\Big\}; A\Big]$$

$$\geq \exp[(\chi(\gamma) - \gamma\chi'(\gamma) - 2\varepsilon)T] \cdot \widehat{P}(A_{1,T} \cap A_{2,T} \cap A)$$

$$\geq \exp[(\chi(\gamma) - \gamma\chi'(\gamma) - 2\varepsilon)T] \cdot \widehat{P}(A_{1,T} \cap A_{2,T} \cap A_{3,T} \cap A_{4,T} \cap A_{5,T})$$

$$\geq \exp[(\chi(\gamma) - \gamma\chi'(\gamma) - 2\varepsilon)T] \cdot (1 - 5\varepsilon).$$

The estimate of $P(A)$ is uniform in $\pi$. Therefore, we see that

$$\Pi(k) \geq \varliminf_{T\to\infty} \frac{1}{T} \log\{\exp[(\chi(\gamma) - \gamma\chi'(\gamma) - 2\varepsilon)T] \cdot (1 - 5\varepsilon)\}$$

$$= \widehat{\chi}(k - \varepsilon) - 2\varepsilon.$$

By continuity of $\widehat{\chi}$ on $(\chi'(-\infty), \chi'(0-))$ and sending $\varepsilon$ to 0, we obtain

$$(3.15) \qquad \Pi(k) \geq \inf_{\gamma<0}\{\chi(\gamma) - \gamma k\}.$$

On the other hand, if $\pi = \widehat{\pi}^{[k]}$, we have

$$P\Big(\frac{\log X_T^{\widehat{\pi}^{[k]}}}{T} \leq k\Big) = P((X_T^{\widehat{\pi}})^{\gamma(k)} \geq e^{\gamma(k)kT})$$

$$\leq E[(X_T^{\widehat{\pi}})^{\gamma(k)}] \cdot e^{-\gamma(k)kT}$$

$$= \exp(\gamma(k)\log x + v(0, y; T; \gamma(k)) - \gamma(k)kT).$$

Therefore, by Proposition 2.3(iii),

$$\varlimsup_{T\to\infty} \frac{1}{T} \log P\Big(\frac{\log X_T^{\widehat{\pi}^{[k]}}}{T} \leq k\Big) \leq \chi(\gamma(k)) - \gamma(k)k$$

$$= \chi(\gamma(k)) - \gamma(k)\chi'(\gamma(k))$$

$$= \inf_{\gamma<0}\{\chi(\gamma) - \gamma k\}.$$

Together with (3.15), we get

$$\Pi(k) = \lim_{T\to\infty} \frac{1}{T} \log P\Big(\frac{\log X_T^{\widehat{\pi}^{[k]}}}{T} \leq k\Big) = \inf_{\gamma<0}\{\chi(\gamma) - \gamma k\}.$$

We have proved (2.48).

We now consider $k < r$. By convexity of $\chi(\cdot)$, we have

$$\chi(-1) \geq \chi(\gamma) + \chi'(\gamma)(-1 - \gamma), \qquad \gamma < -1.$$

That is,

$$\chi(\gamma) - \gamma k \leq \chi(-1) + \chi'(\gamma) + \gamma(\chi'(\gamma) - k).$$

Since $\chi'(\gamma)$ is bounded, $\chi'(\gamma) \to r$ as $\gamma \to -\infty$ (Proposition 2.7), we see

$$\chi(\gamma) - \gamma k \to -\infty, \qquad \gamma \to -\infty.$$

In particular,

$$\inf_{\gamma < 0} \{\chi(\gamma) - \gamma k\} = -\infty.$$

On the other hand, we take $\pi = 0$. Then

$$X_T^0 = x \exp(rT),$$

$$P\left(\frac{\log X_T^0}{T} \le k\right) = 0 \qquad \text{if } T \text{ is sufficiently large.}$$

In particular, $\Pi(k) = -\infty$.

We now assume that $B$ is stable. Consider $\frac{d}{dt} e^{tK_1^*} \overline{P} e^{tK_1}$ and use (2.34). We can show

$$\overline{P} - e^{tK_1^*} \overline{P} e^{tK_1} = \int_0^t e^{sK_1^*} (\overline{P} \Lambda N^{-1} \Lambda^* \overline{P} - C^*C) e^{sK_1} \, ds.$$

Since $K_1$ is stable if $\gamma < 0$ is near 0 [see (2.28)], when $t$ tends to infinite, we have

$$-\overline{P} + \int_0^\infty e^{sK_1^*} \overline{P} \Lambda N^{-1} \Lambda^* \overline{P} e^{sK_1} \, ds = \int_0^\infty e^{sK_1^*} C^*C e^{sK_1} \, ds.$$

We want to let $\gamma \to 0-$ for $\overline{P} = \overline{P}(\gamma)$. By Proposition 2.3(i), $-\overline{P}(\gamma)$ is non-negative. The above relation implies that $-\overline{P}(\gamma)$ is bounded above as $\gamma \to 0$. By (4.9) in the next section, $\overline{P}(\gamma)$ is nondecreasing. Therefore,

$$\overline{P}(0-) = \lim_{\gamma \to 0-} \overline{P}(\gamma)$$

exists and is nonpositive definite. Since $C(\gamma) \to 0$ as $\gamma \to 0-$, we take the limit in the above relation and we easily see $\overline{P}(0-) = 0$.

By (2.35),

$$B^* \overline{q}(0-) = 0.$$

We have $\overline{q}(0-) = 0$. From (2.37), we have $\chi(0-) = 0$. We have proved (2.50).

Using (4.9), we obtain

$$\frac{d\overline{P}}{d\gamma}(0-) = \int_0^\infty e^{sB^*} A^* (\Sigma\Sigma^*)^{-1} A e^{sB} \, ds.$$

By (4.12), we have

$$\frac{d\overline{q}}{d\gamma}(0-) = -(B^*)^{-1} \left[ \frac{d\overline{P}}{d\gamma}(0-) b + A^* (\Sigma\Sigma^*)^{-1} (a - r\mathbf{1}) \right].$$

By (4.13), we have

$$\chi'(0-) = \frac{1}{2}\operatorname{tr}\left(\Lambda\Lambda^*\frac{d\overline{P}}{d\gamma}(0-)\right) + \frac{d\overline{q}}{d\gamma}(0-)^*b + \frac{1}{2}(a-r\mathbf{1})^*(\Sigma\Sigma^*)^{-1}(a-r\mathbf{1}) + r.$$

Using the above expression for $\frac{d\overline{P}}{d\gamma}(0-)$, we can show

$$B^*\frac{d\overline{P}}{d\gamma}(0-) + \frac{d\overline{P}}{d\gamma}(0-)B = -A^*(\Sigma\Sigma^*)^{-1}A.$$

Multiplying this by $(B^*)^{-1}$ on the left and $B^{-1}$ on the right, we have

$$(B^*)^{-1}\frac{d\overline{P}}{d\gamma}(0-) + \frac{d\overline{P}}{d\gamma}(0-)B^{-1} = -(B^*)^{-1}A^*(\Sigma\Sigma^*)^{-1}AB^{-1}.$$

Then

$$(B^*)^{-1}\frac{d\overline{P}}{d\gamma}(0-)b \cdot b$$

$$= -\frac{1}{2}(AB^{-1}b)^*(\Sigma\Sigma^*)^{-1}AB^{-1}b,$$

$$\frac{d\overline{q}}{d\gamma}(0-)^*b + \frac{1}{2}(a-r\mathbf{1})^*(\Sigma\Sigma^*)^{-1}(a-r\mathbf{1})$$

$$= \frac{1}{2}(AB^{-1}b - (a-r\mathbf{1}))^*(\Sigma\Sigma^*)^{-1}(AB^{-1}b - (a-r\mathbf{1})).$$

An expression of $\chi'(0-)$ as in (2.51) follows.

Finally, we show

$$\inf_{\gamma<0}\{\chi(\gamma) - \gamma k\} = 0, \qquad k \geq \chi'(0-).$$

Let $k \geq \chi'(0-)$. Since $\chi(0-) = 0$, we have

$$\inf_{\gamma<0}\{\chi(\gamma) - \gamma k\} \leq 0.$$

Take $k_1 < \chi'(0-)$. Then

$$\Pi(k_1) = \inf_{\gamma<0}\{\chi(\gamma) - \gamma k_1\} \leq \inf_{\gamma<0}\{\chi(\gamma) - \gamma k\},$$

$$\Pi(k_1) = \chi(\gamma(k_1)) - \gamma(k_1)k_1,$$

where

$$\chi'(\gamma(k_1)) = k_1.$$

$\gamma(k_1) \to 0-$ as $k_1 \to \chi'(0-)$. Then $\Pi(k_1) \to 0$ as $k_1 \to \chi'(0-)$. Therefore,

$$\inf_{\gamma<0}\{\chi(\gamma) - \gamma k\} \geq 0;$$

hence

$$\inf_{\gamma<0}\{\chi(\gamma)-\gamma k\}=0.$$

This completes the proof of our main theorem. The following is a lemma used in the above proof.

LEMMA 3.1. *Assume (A1) and (A2). Then, for all $p \geq 1$, there is $M = M(\gamma,p)$ such that*

(3.16)                    $$\widehat{E}[|Y_t|^{2p}] \leq |y|^{2p} + M, \qquad t \geq 0.$$

PROOF. Let $K_\gamma$ be the positive definite matrix in Proposition 2.6. We apply Itô's differential rule to $(Y_t^* K_\gamma Y_t)^p$. $Y_t$ is given in (2.45) and we consider the probability $\widehat{P}$ [in (2.43)]. To simplify the calculation, we assume $K_\gamma = I$ (identity matrix) in the following. Noting that (2.46) holds for $K_\gamma = I$, we have

$$\begin{aligned}
d|Y_t|^{2p} &= 2p|Y_t|^{2(p-1)}Y_t^*\Lambda \, d\widehat{W}_t + 2p|Y_t|^{2(p-1)}Y_t^*\beta_\gamma(Y_t)\, dt \\
&\quad + p|Y_t|^{2(p-1)}\operatorname{tr}\Lambda\Lambda^*\, ds + 2p(p-1)|Y_t|^{2(p-2)}Y_t^*\Lambda\Lambda^* Y_t\, dt \\
&\leq 2p|Y_t|^{2(p-1)}Y_t^*\Lambda \, d\widehat{W}_t \\
&\quad + p|Y_t|^{2(p-1)}\{-2c_1(\gamma)|Y_t|^2 + 2c_2(\gamma)+c_3\}\, dt
\end{aligned}$$

for some $c_3 > 0$. Moreover, if we take $\theta$ as $0 < \theta < 2pc_1(\gamma)$, then we have

$$\begin{aligned}
d|Y_t|^{2p}e^{\theta t} &\leq 2pe^{\theta t}|Y_t|^{2(p-1)}Y_t^*\Lambda \, d\widehat{W}_t \\
&\quad + e^{\theta t}|Y_t|^{2(p-1)}[p\{-2c_1(\gamma)|Y_t|^2+2c_2(\gamma)+c_3\}+\theta|Y_t|^2]\, dt \\
&\leq 2pe^{\theta t}|Y_t|^{2(p-1)}Y_t^*\Lambda \, d\widehat{W}_t + e^{\theta t}\overline{M}\, dt,
\end{aligned}$$

where

$$\overline{M} := \sup_y |y|^{2(p-1)}[p\{-2c_1(\gamma)|y|^2+2c_2(\gamma)+c_3\}+\theta|y|^2].$$

Therefore, we obtain

$$|Y_t|^{2p} \leq |y|^{2p} + \frac{\overline{M}}{\theta}(1-e^{-\theta t}) + 2p\int_0^t e^{-\theta(t-s)}|Y_s|^{2(p-1)}Y_s^*\Lambda \, d\widehat{W}_s.$$

From this, by an argument using a stopping time, we obtain (3.16) with $M = \overline{M}/\theta$. □

REMARK 3.1. The proof of Theorem 2.1 is based on the results in Propositions 2.1–2.7. It is interesting to verify these results for more general non-Gaussian models. We make the following observation. If we assume:

(i) $\alpha, \sigma, \beta, \lambda$ are globally Lipschitz;

(ii) $\mu_1|\xi|^2 \leq \xi^*\sigma\sigma^*(y)\xi \leq \mu_2|\xi|^2, \mu_1, \mu_2 > 0$;

(iii) $\nu_1|\xi|^2 \leq \xi^*\lambda\lambda^*(y)\xi \leq \nu_2|\xi|^2, \nu_1, \nu_2 > 0$; then there is $\chi^*(\gamma)$ such that (2.23) has a solution for $\chi \geq \chi^*(\gamma)$ (see Theorem 2.6 in [22]). Moreover, under certain conditions, there is a unique solution $\xi_\gamma(y)$ with $\xi_\gamma(0) = 0$ for $\chi = \chi^*(\gamma)$ (see Theorem 3.8 in [22]). Assume $\Lambda\Lambda^*$ is positive, and for Gaussian cases studied here that satisfy the conditions (A1), (A2), (A3), $\chi^*(\gamma) = \chi(\gamma), \xi_\gamma^*(y) = \xi_\gamma(y)$, which are defined in Proposition 2.3. It is an interesting problem to find conditions on coefficients that prove other propositions.

## 4. Proof of propositions.

4.1. *Finite time horizon problem.* In this subsection, we prove Propositions 2.1 and 2.2. We follow closely the arguments of Kuroda and Nagai [25].

First of all, we attack the finite time horizon problem (2.9). Then the solution $v$ of the Bellman equation (2.21) can be expressed as quadratic function such that

$$v(t, y) = \tfrac{1}{2}y^*P(t)y + q(t)^*y + h(t),$$

provided that equation (2.27) has a solution. Here $q$ and $h$ are solutions of (2.29) and (2.30), respectively. We recall the following result (5.2) in [13], Theorem IV. If $\gamma < 0$, then we see that (2.27) has the unique solution $P(t) \leq 0$. See also Remark 1 of Section 2 in [25], but notice that $\gamma$ and the solution $P(t)$ of (2.27) correspond to $-\frac{\theta}{2}$ and $-\frac{\theta}{2}P(t)$ in [25], respectively. Therefore, as in the proof of Theorem 2.1 in [25], we obtain Propositions 2.1 and 2.2.

4.2. *Asymptotics as $T \to \infty$.* In this section we prove Propositions 2.3 and 2.4. Here we recall the following theorem (see Theorem 4.1 and Lemma 5.2 in [38]).

THEOREM 4.1 [38]. *Assume that $N > 0$ and $(K_1, \Lambda)$ is stabilizable, then for the solution $Q(t) = Q(t; T)$ of*

$$(4.1) \qquad \begin{cases} \dot{Q}(t) + K_1^*Q(t) + Q(t)K_1 - Q(t)\Lambda N^{-1}\Lambda^*Q(t) + C^*C = 0, \\ Q(T) = 0. \end{cases}$$

$\exists \lim_{T \to \infty} Q(t; T) = \overline{Q}$, *and $\overline{Q}$ satisfies the algebraic Riccati equation,*

$$(4.2) \qquad K_1^*\overline{Q} + \overline{Q}K_1 - \overline{Q}\Lambda N^{-1}\Lambda^*\overline{Q} + C^*C = 0.$$

*Moreover, if $(C, K_1)$ is detectable, then $K_1^* - \Lambda N^{-1}\Lambda^*\overline{Q}$ is stable and the nonnegative definite solution of (4.2) is unique.*

REMARK 4.1. The pair $(L, M)$ of the $n \times n$ matrix $L$ and the $n \times l$ matrix $M$ is stabilizable if there exists a $l \times n$ matrix $K$ such that $L - MK$ is stable. The pair $(L, F)$ of the $l \times n$ matrix $L$ and the $n \times n$ matrix $F$ is called detectable if $(F^*, L^*)$ is stabilizable.

PROOF OF PROPOSITION 2.3(i) AND (ii). Let us set $Q(t) = -P(t)$. Then we see that $Q$ satisfies (4.1). Let $K_1, C$ be given in (2.28). Take $K = \Sigma^*(\Sigma\Sigma^*)^{-1}A$. Then $K_1 - \Lambda K = B - \Lambda\Sigma^*(\Sigma\Sigma^*)^{-1}A = G$ is stable. Note also that if we set $K_2 := \sqrt{\frac{-1}{\gamma(1-\gamma)}}\Lambda^*$, then $K_1^* - CK_2 = G^*$ is stable [$G$ is given in (A2)]. Therefore, we see that $(K_1, \Lambda)$ is stabilizable, $(C, K_1)$ is detectable and Theorem 4.1 applies. Then we can follow the arguments of Proposition 2.2 in Kuroda and Nagai [25]; we omit the details here.

PROOF OF PROPOSITION 2.3(iii). Note that $P(t; T)$ is uniformly bounded with respect to $t$ and $T$ (see Remark 1 of Section 4 in [25]). Moreover, since $K_1 + \Lambda N^{-1}\Lambda^* P(t; T)$ converges to a stable matrix $K_1 + \Lambda N^{-1}\Lambda^*\overline{P}$ as $T \to \infty$, we can show by (2.29) that $q(t; T)$ is uniformly bounded with respect to $t$ and $T$. See Lemma 4.4 in [25]. Therefore we obtain (2.38).  □

PROOF OF PROPOSITION 2.4. Let us assume (A1)–(A3). Then we have the solution $\overline{P}$ of (2.34) and $\overline{q}$ of (2.35). Moreover, the pair $(\xi(y), \chi(\gamma))$ of $\xi(y)$, defined by (2.33) and $\chi(\gamma)$ defined by (2.37), satisfies (2.23) (cf. Section 2 in Kuroda and Nagai [25]). Here we note that from (A2), $K_1 + \Lambda N^{-1}\Lambda^*\overline{P}$ is stable (by Theorem 4.1). Moreover, from (A3) we see that the variance of $Y_t$ under $\widehat{P}$ [defined in (2.43)] is nondegenerate (see Lemma 5.1 of Kuroda and Nagai [25]). Therefore, we see that $Y_t$ given by (2.45) is $\widehat{P}$-ergodic. The rest of the proof follows closely the arguments of Theorem 3.8 in Kaise and Sheu [22] and is omitted here.  □

4.3. *Differentiability with respect to HARA parameter $\gamma$.* In this subsection, we shall prove Proposition 2.5(i). Let us first note that the Riccati differential equation (2.1) can be solved by considering a Hamiltonian system. Indeed, introduce a Hamiltonian matrix $\mathcal{H}$ defined by

$$(4.3) \qquad \mathcal{H} = \begin{pmatrix} -K_1 & \Lambda N^{-1}\Lambda^* \\ C^*C & K_1^* \end{pmatrix},$$

and consider the ordinary differential equation

$$(4.4) \quad \begin{pmatrix} \dot{U}(t) \\ \dot{V}(t) \end{pmatrix} = \mathcal{H}\begin{pmatrix} \dot{U}(t) \\ \dot{V}(t) \end{pmatrix}, \qquad 0 \leq t \leq T, \qquad \begin{pmatrix} \dot{U}(0) \\ \dot{V}(0) \end{pmatrix} = \begin{pmatrix} E_n \\ 0 \end{pmatrix}.$$

See Chapter V in [1]. Note that $U$ and $V$ are $n \times n$ matrix valued functions on $0 \leq t \leq T$, and $E_n$ is the $n \times n$ unit matrix. Then it is known that $U(t)$ is invertible, and $W(t) := V(t)U(t)^{-1}$ is the solution to the Riccati differential equation

$$(4.5) \quad \begin{cases} \dot{W}(t) = K_1^* W(t) + W(t)K_1 - W(t)\Lambda N^{-1}\Lambda^* W(t) + C^*C, \\ W(0) = 0. \end{cases}$$

Then we see that, by setting $\widehat{P}(t) := -P(T-t; T)$, we have $\widehat{P}(t) = W(t)$.

LEMMA 4.1. *The solution* $P(t) := P(t; T; \gamma)$ *to the Riccati equation (2.27) is in* $\mathbf{C}^1$*-class with respect to* $\gamma$.

PROOF. The Hamiltonian matrix $\mathcal{H}$ defined by (4.3) is smooth with respect to $\gamma$ and so is the solution $(U(t), V(t))^*$ of (4.4). Moreover, $U(t)$ is invertible. Therefore, $U(t)^{-1}$ is in $\mathbf{C}^1$-class and

$$\frac{\partial U(t)^{-1}}{\partial \gamma} = -U(t)^{-1} \frac{\partial U(t)}{\partial \gamma} U(t)^{-1}.$$

Thus we see that $W(t) = V(t)U(t)^{-1}$ is in $\mathbf{C}^1$-class with respect to $\gamma$. Hence we conclude our present lemma. $\square$

Now, let us rewrite (2.27) as

$$(4.6) \quad \begin{aligned} &\dot{P}(t) + (K_1 + \Lambda N^{-1}\Lambda^* P(t))^* P(t) + P(t)(K_1 + \Lambda N^{-1}\Lambda^* P(t)) \\ &\quad - P(t)\Lambda N^{-1}\Lambda^* P(t) - C^*C = 0. \end{aligned}$$

Then by differentiating (4.6) with respect to $\gamma$, we obtain

$$(4.7) \quad \begin{aligned} &\frac{d}{dt}\left(\frac{\partial P}{\partial \gamma}\right) + (K_1 + \Lambda N^{-1}\Lambda^* P(t))^* \frac{\partial P}{\partial \gamma} \\ &\quad + \frac{\partial P}{\partial \gamma}(K_1 + \Lambda N^{-1}\Lambda^* P(t)) \\ &\quad + \frac{1}{(1-\gamma)^2}(\Sigma\Lambda^* P(t) + A)^*(\Sigma\Sigma^*)^{-1}(\Sigma\Lambda^* P(t) + A) = 0. \end{aligned}$$

Then we obtain the following lemma.

LEMMA 4.2. *Assume (A1) and (A2). Then the solution* $\frac{\partial P}{\partial \gamma}(t; T; \gamma)$ *of (4.7) converges to* $\frac{d\overline{P}}{d\gamma}$ *which satisfies*

$$(4.8) \quad \begin{aligned} &(K_1 + \Lambda N^{-1}\Lambda^*\overline{P})^* \frac{d\overline{P}}{d\gamma} + \frac{d\overline{P}}{d\gamma}(K_1 + \Lambda N^{-1}\Lambda^*\overline{P}) \\ &\quad + \frac{1}{(1-\gamma)^2}(\Sigma\Lambda^*\overline{P} + A)^*(\Sigma\Sigma^*)^{-1}(\Sigma\Lambda^*\overline{P} + A) = 0. \end{aligned}$$

*Moreover, we obtain the following expression:*

$$(4.9) \quad \frac{d\overline{P}}{d\gamma} = \frac{1}{(1-\gamma)^2} \int_0^\infty e^{s(K_1 + \Lambda N^{-1}\Lambda^* \overline{P})^*} (\Sigma\Lambda^* \overline{P} + A)^*$$
$$\times (\Sigma\Sigma^*)^{-1} (\Sigma\Lambda^* \overline{P} + A) e^{s(K_1 + \Lambda N^{-1}\Lambda^* \overline{P})} \, ds.$$

PROOF. Note that $K_1 + \Lambda N^{-1}\Lambda^* P(t; T)$ converges to the stable matrix $K_1 + \Lambda N^{-1}\Lambda^* \overline{P}$. We can see that similar to Lemma 4.4 of [25], $\frac{\partial P}{\partial \gamma}(t; T)$ converges to a matrix $\left(\overline{\frac{\partial P}{\partial \gamma}}\right)$ which satisfies

$$(K_1 + \Lambda N^{-1}\Lambda^* \overline{P})^* \left(\overline{\frac{\partial P}{\partial \gamma}}\right)$$

$$+ \left(\overline{\frac{\partial P}{\partial \gamma}}\right)(K_1 + \Lambda N^{-1}\Lambda^* \overline{P})$$

$$+ \frac{1}{(1-\gamma)^2} (\Sigma\Lambda^* \overline{P} + A)^* (\Sigma\Sigma^*)^{-1} (\Sigma\Lambda^* \overline{P} + A) = 0.$$

Then

$$\frac{d}{dt} \left\{ e^{t(K_1 + \Lambda N^{-1}\Lambda^* \overline{P})^*} \left(\overline{\frac{\partial P}{\partial \gamma}}\right) e^{t(K_1 + \Lambda N^{-1}\Lambda^* \overline{P})} \right\}$$

$$= e^{t(K_1 + \Lambda N^{-1}\Lambda^* \overline{P})^*} \left( (K_1 + \Lambda N^{-1}\Lambda^* \overline{P})^* \left(\overline{\frac{\partial P}{\partial \gamma}}\right) \right.$$

$$\left. + \left(\overline{\frac{\partial P}{\partial \gamma}}\right)(K_1 + \Lambda N^{-1}\Lambda^* \overline{P}) \right) e^{t(K_1 + \Lambda N^{-1}\Lambda^* \overline{P})}$$

$$= -\frac{1}{(1-\gamma)^2} e^{t(K_1 + \Lambda N^{-1}\Lambda^* \overline{P})^*} (\Sigma\Lambda^* \overline{P} + A)^* (\Sigma\Sigma^*)^{-1}$$

$$\times (\Sigma\Lambda^* \overline{P} + A) e^{t(K_1 + \Lambda N^{-1}\Lambda^* \overline{P})}.$$

Integrating over $t$, then $\left(\overline{\frac{\partial P}{\partial \gamma}}\right)$ satisfies

$$\left(\overline{\frac{\partial P}{\partial \gamma}}\right) - e^{t(K_1 + \Lambda N^{-1}\Lambda^* \overline{P})^*} \left(\overline{\frac{\partial P}{\partial \gamma}}\right) e^{t(K_1 + \Lambda N^{-1}\Lambda^* \overline{P})}$$

$$= \frac{1}{(1-\gamma)^2} \int_0^t e^{s(K_1 + \Lambda N^{-1}\Lambda^* \overline{P})^*} (\Sigma\Lambda^* \overline{P} + A)^*$$

$$\times (\Sigma\Sigma^*)^{-1} (\Sigma\Lambda^* \overline{P} + A) e^{s(K_1 + \Lambda N^{-1}\Lambda^* \overline{P})} \, ds.$$

We see that

$$\left(\overline{\frac{\partial P}{\partial \gamma}}\right) = \frac{1}{(1-\gamma)^2} \int_0^\infty e^{s(K_1 + \Lambda N^{-1}\Lambda^* \overline{P})^*} (\Sigma\Lambda^* \overline{P} + A)^*$$

(4.10)
$$\times (\Sigma\Sigma^*)^{-1}(\Sigma\Lambda^*\overline{P} + A)e^{s(K_1 + \Lambda N^{-1}\Lambda^*\overline{P})}\,ds.$$

On the other hand, we have

$$\overline{P}(\gamma + \Delta) - \overline{P}(\gamma) = \lim_{T\to\infty}\{P(t;T;\gamma + \Delta) - P(t;T;\gamma)\}$$

$$= \lim_{T\to\infty}\int_\gamma^{\gamma+\Delta}\frac{\partial P}{\partial\gamma}(t;T;u)\,du$$

$$= \int_\gamma^{\gamma+\Delta}\left(\overline{\frac{\partial P}{\partial\gamma}}\right)(u)\,du.$$

From (4.10), $(\overline{\frac{\partial P}{\partial\gamma}})$ is continuous with respect to $\gamma$. Therefore, we see that $\overline{P}$ is differentiable with respect to $\gamma$, and $\frac{d\overline{P}}{d\gamma}(\gamma) = (\overline{\frac{\partial P}{\partial\gamma}})(\gamma)$. $\quad\square$

As for differentiability of $\overline{q}$ with respect to $\gamma$, we can see this is similar to Lemma 4.2. Indeed, (2.29) is a linear equation and its coefficients are all in $\mathbf{C}^1$-class with respect to $\gamma$. Therefore, the solution $q(t)$ of (2.29) is in $\mathbf{C}^1$-class with respect to $\gamma$, and we have

$$\frac{d}{dt}\left(\frac{\partial q}{\partial\gamma}\right) + (K_1 + \Lambda N^{-1}\Lambda^* P(t))^*\frac{\partial q}{\partial\gamma} + \frac{\partial(K_1 + \Lambda N^{-1}\Lambda^* P(t))^*}{\partial\gamma}q(t)$$

(4.11)
$$+ \frac{\partial P}{\partial\gamma}b + \frac{1}{(1-\gamma)^2}(\Sigma\Lambda^* P(t) + A)^*(\Sigma\Sigma^*)^{-1}(a - r\mathbf{1})$$

$$+ \frac{\gamma}{1-\gamma}\frac{\partial P}{\partial\gamma}\Lambda\Sigma^*(\Sigma\Sigma^*)^{-1}(a - r\mathbf{1}) = 0.$$

Thus we have the following lemma, similar to Lemma 4.2.

LEMMA 4.3. *Under assumptions (A1) and (A2), as $T\to\infty$, $\frac{\partial q}{\partial\gamma}(t;T;\gamma)$, the solution of (4.11) converges to $\frac{d\overline{q}}{d\gamma}$ which satisfies*

(4.12)
$$(K_1 + \Lambda N^{-1}\Lambda^*\overline{P})^*\frac{d\overline{q}}{d\gamma} + \Phi\left(\gamma, \frac{d\overline{P}}{d\gamma}, \overline{P}, \overline{q}\right) = 0,$$

*where*

$$\Phi\left(\gamma, \frac{d\overline{P}}{d\gamma}, \overline{P}, \overline{q}\right) := \left[\frac{1}{(1-\gamma)^2}\Lambda\Sigma^*(\Sigma\Sigma^*)^{-1}(A + \Sigma\Lambda^*\overline{P}) + \Lambda N^{-1}\Lambda^*\frac{d\overline{P}}{d\gamma}\right]^*\overline{q}$$

$$+ \frac{d\overline{P}}{d\gamma}b + \frac{1}{(1-\gamma)^2}(A + \Sigma\Lambda^*\overline{P})^*(\Sigma\Sigma^*)^{-1}(a - r\mathbf{1})$$

$$+ \frac{\gamma}{1-\gamma}\frac{d\overline{P}}{d\gamma}\Lambda\Sigma^*(\Sigma\Sigma^*)^{-1}(a - r\mathbf{1}).$$

Differentiability of $\chi(\gamma)$ is directly seen from (2.37). Indeed, we have

$$
\begin{aligned}
\frac{d\chi}{d\gamma} = \frac{1}{2}\operatorname{tr}&\left(\Lambda\Lambda^*\frac{d\overline{P}}{d\gamma}\right) + \overline{q}^*\Lambda\Lambda^*\frac{d\overline{q}}{d\gamma} + b^*\frac{d\overline{q}}{d\gamma} + r \\
(4.13)\qquad & + \frac{1}{2(1-\gamma)^2}(a - r\mathbf{1} + \Sigma\Lambda^*\overline{q})^*(\Sigma\Sigma^*)^{-1}(a - r\mathbf{1} + \Sigma\Lambda^*\overline{q}) \\
& + \frac{\gamma}{1-\gamma}(a - r\mathbf{1} + \Sigma\Lambda^*\overline{q})^*(\Sigma\Sigma^*)^{-1}\Sigma\Lambda^*\frac{d\overline{q}}{d\gamma}.
\end{aligned}
$$

The following lemma is a direct consequence of (4.9), (4.12) and (4.13).

LEMMA 4.4. *Under assumptions (A1) and (A2), $\frac{d\overline{P}}{d\gamma}, \frac{d\overline{q}}{d\gamma}$ and $\frac{d\chi}{d\gamma}$ are differentiable with respect to $\gamma$.*

PROOF. Differentiability of $\frac{d\overline{P}}{d\gamma}$ is seen by looking at (4.9). As for $\frac{d\overline{q}}{d\gamma}$, from (4.12) we obtain

$$
\frac{d\overline{q}}{d\gamma} = -[(K_1 + \Lambda N^{-1}\Lambda^*\overline{P})^*]^{-1}\Phi\left(\gamma, \frac{d\overline{P}}{d\gamma}, \overline{P}, \overline{q}\right),
$$

and so it turns out to be differentiable. From these facts and (4.13), differentiability of $\frac{d\chi}{d\gamma}$ follows.   □

4.4. *Convexity of $\chi$.*   In this subsection, we shall show Proposition 2.5(ii).

PROOF OF PROPOSITION 2.5(ii).   Note that $K_1 + \Lambda N^{-1}\Lambda^*\overline{P}$ is stable under assumption (A2). In the previous subsection, namely the proof of Proposition 2.5(i), we have shown in Lemma 4.4 that under assumptions (A1) and (A2), $\overline{P}, \overline{q}$ and $\chi$ are twice differentiable with respect to $\gamma$ and so is $\xi$. Recall that

$$
\frac{d\chi}{d\gamma}(\gamma) = \frac{1}{2}\operatorname{tr}(\Lambda\Lambda^* D^2\eta(y)) + \beta_\gamma(y)^* D\eta(y) + V_1,
$$

where $\eta := \frac{\partial\xi}{\partial\gamma}$, $\beta_\gamma(y)$ is given by (2.40) and $V_1$ is defined by (3.5). Moreover, setting $\zeta := \frac{\partial^2\xi}{\partial\gamma^2}$, we have

$$
\begin{aligned}
\frac{d^2\chi}{d\gamma^2}(\gamma) = \frac{1}{2}\operatorname{tr}&(\Lambda\Lambda^* D^2\zeta(y)) + \beta_\gamma(y)^* D\zeta(y) + \left(\frac{\partial\beta_\gamma(y)}{\partial\gamma}\right)^* D\eta(y) \\
& + \frac{1}{(1-\gamma)^3}\{\Sigma\Lambda^* D\xi(y) + (a + Ay - r\mathbf{1})\}^* \\
& \times (\Sigma\Sigma^*)^{-1}\{\Sigma\Lambda^* D\xi(y) + (a + Ay - r\mathbf{1})\}
\end{aligned}
$$

$$+ \frac{1}{(1-\gamma)^2}(D\eta)^*(y)\Lambda\Sigma^*(\Sigma\Sigma^*)^{-1}$$
$$\times \{\Sigma\Lambda^* D\xi(y) + (a + Ay - r\mathbf{1})\}.$$

Using

$$\frac{\partial \beta_\gamma}{\partial \gamma}(y) = \frac{1}{(1-\gamma)^2}\Lambda\Sigma^*(\Sigma\Sigma^*)^{-1}\{\Sigma\Lambda^* D\xi(y) + (a + Ay - r\mathbf{1})\}$$
$$+ \Lambda N^{-1}\Lambda^* D\eta(y),$$

we obtain

$$\frac{d^2\chi}{d\gamma^2}(\gamma) = \frac{1}{2}\operatorname{tr}(\Lambda\Lambda^* D^2\zeta(y)) + \beta_\gamma(y)^* D\zeta + \ell(y;\gamma),$$

where

$$\ell(y;\gamma) := (D\eta)^*(y)\Lambda(I - \Sigma^*(\Sigma\Sigma^*)^{-1}\Sigma)\Lambda^* D\eta(y)$$
$$+ \frac{1}{(1-\gamma)^3}\{(1-\gamma)\Sigma\Lambda^* D\eta(y) + \Sigma\Lambda^* D\xi(y) + (a + Ay - r\mathbf{1})\}^*$$
$$\times (\Sigma\Sigma^*)^{-1}\{(1-\gamma)\Sigma\Lambda^* D\eta(y) + \Sigma\Lambda^* D\xi(y) + (a + Ay - r\mathbf{1})\}.$$

Note that $|\ell(y;\gamma)| \leq K(1 + |y|^2)$. From (3.16) we can see that

$$\frac{d^2\chi}{d\gamma^2}(\gamma) = \lim_{T\to\infty}\frac{1}{T}\widehat{E}\left[\int_0^T \ell(Y_s;\gamma)\,ds\right] < \infty.$$

Since $\ell(y;\gamma) \geq 0$, we conclude $\frac{d^2\chi}{d\gamma^2}(\gamma) \geq 0$. Therefore, Proposition 2.5(ii) is obtained. $\square$

4.5. *Proof of Proposition 2.6.* Recall that $K_1 + \Lambda N^{-1}\Lambda^*\overline{P}$ is stable from Theorem 4.1. Let us set $\overline{G} := K_1 + \Lambda N^{-1}\Lambda^*\overline{P}$, and consider

$$K := \int_0^\infty e^{t\overline{G}^*}e^{t\overline{G}}\,dt > 0.$$

Then $K$ satisfies the following equation:

$$\overline{G}^* K + K\overline{G} = -I.$$

Therefore, we have

$$\langle K\overline{G}y, y\rangle + \langle y, K\overline{G}y\rangle = -\langle y, y\rangle,$$
$$\langle \overline{G}y, Ky\rangle = -\tfrac{1}{2}\langle y, y\rangle.$$

Since

$$\beta_\gamma(y) = \overline{G}y + f_\gamma$$

[see (2.40)], we can deduce (2.46) after some calculation. Here $K_\gamma = K$, $c_1(\gamma) = 1/4$, and $c_2(\gamma) = |Kf_\gamma|^2$.

4.6. *Asymptotics as $\gamma \to -\infty$.*   In this subsection we shall consider asymptotic behavior of $\frac{d\chi}{d\gamma}(\gamma)$ as $\gamma \to -\infty$, and obtain Proposition 2.7.

PROOF OF PROPOSITION 2.7.   We first consider $\overline{P}(\gamma)$ which is a solution of the algebraic Riccati equation (2.34). We then note that $\frac{d\overline{P}}{d\gamma} \geq 0$ holds, and $\overline{P}$ is bounded because of (4.9) and (2.36). We set $\overline{P}_{(-\infty)} := \lim_{\gamma \to -\infty} \overline{P}(\gamma)$. Now we rewrite (2.34) as

$$(4.14) \quad \begin{aligned} (K_1 - \Lambda K)^* \overline{P} + \overline{P}(K_1 - \Lambda K) + (\Lambda^* \overline{P} + NK)^* N^{-1}(\Lambda^* \overline{P} + NK) \\ - K^* NK - C^* C = 0, \end{aligned}$$

where

$$(4.15) \qquad K := \frac{1}{1-\gamma} \Sigma^* (\Sigma \Sigma^*)^{-1} A.$$

Noting that

$$N = I - \gamma \Sigma^* (\Sigma \Sigma^*)^{-1} \Sigma,$$

$$\lim_{\gamma \to -\infty} K = 0, \qquad \lim_{\gamma \to -\infty} K_1 = G,$$

$$\lim_{\gamma \to -\infty} N^{-1} = I - \Sigma^* (\Sigma \Sigma^*)^{-1} \Sigma := \widehat{N}_{(-\infty)},$$

$$\lim_{\gamma \to -\infty} (K^* NK + C^* C) = A^* (\Sigma \Sigma^*)^{-1} A,$$

where $G = B - \Lambda \Sigma^* (\Sigma \Sigma^*)^{-1} A$ [see (A2)]. We obtain

$$(4.16) \quad \begin{aligned} G^* \overline{P}_{(-\infty)} + \overline{P}_{(-\infty)} G \\ + \overline{P}_{(-\infty)} \Lambda \widehat{N}_{(-\infty)} \Lambda^* \overline{P}_{(-\infty)} - A^* (\Sigma \Sigma^*)^{-1} A = 0. \end{aligned}$$

Moreover, we rewrite (4.16) as

$$\begin{aligned} (G + \Lambda \widehat{N}_{(-\infty)} \Lambda^* \overline{P}_{(-\infty)})^* \overline{P}_{(-\infty)} + \overline{P}_{(-\infty)} (G + \Lambda \widehat{N}_{(-\infty)} \Lambda^* \overline{P}_{(-\infty)}) \\ - (\overline{P}_{(-\infty)} \Lambda \widehat{N}_{(-\infty)} \Lambda^* \overline{P}_{(-\infty)} + A^* (\Sigma \Sigma^*)^{-1} A) = 0. \end{aligned}$$

We consider

$$\frac{d}{dt} e^{t(G + \Lambda \widehat{N}_{(-\infty)} \Lambda^* \overline{P}_{(-\infty)})^*} \overline{P}_{(-\infty)} e^{t(G + \Lambda \widehat{N}_{(-\infty)} \Lambda^* \overline{P}_{(-\infty)})}$$

and using the above relation, we can show

$$(4.17) \quad \begin{aligned} \int_0^\infty e^{s(G + \Lambda \widehat{N}_{(-\infty)} \Lambda^* \overline{P}_{(-\infty)})^*} \overline{P}_{(-\infty)} \\ \times \Lambda \widehat{N}_{(-\infty)} \Lambda^* \overline{P}_{(-\infty)} e^{s(G + \Lambda \widehat{N}_{(-\infty)} \Lambda^* \overline{P}_{(-\infty)})} \, ds \\ \leq -\overline{P}_{(-\infty)}. \end{aligned}$$

Here we use $\overline{P}_{(-\infty)} \le 0$. Since $G^*$ is stable,

$$(G^* + (\Lambda \widehat{N}_{(-\infty)} \Lambda^* \overline{P}_{(-\infty)})^*, (\widehat{N}_{(-\infty)} \Lambda^* \overline{P}_{(-\infty)})^*)$$

is stabilizable which means that $(\widehat{N}_{(-\infty)} \Lambda^* \overline{P}_{(-\infty)}, G + \Lambda \widehat{N}_{(-\infty)} \Lambda^* \overline{P}_{(-\infty)})$ is detectable. Therefore, noting that $(\widehat{N}_{(-\infty)})^2 = \widehat{N}_{(-\infty)}$ and that

$$\left\| \int_0^\infty e^{s(G + \Lambda \widehat{N}_{(-\infty)} \Lambda^* \overline{P}_{(-\infty)})^*} \overline{P}_{(-\infty)} \Lambda \widehat{N}_{(-\infty)} \Lambda^* \overline{P}_{(-\infty)} e^{s(G + \Lambda \widehat{N}_{(-\infty)} \Lambda^* \overline{P}_{(-\infty)})} \, ds \right\|$$

is bounded because of (4.17), we see that $G + \Lambda \widehat{N}_{(-\infty)} \Lambda^* \overline{P}_{(-\infty)}$ is stable (see [39], Proposition 3.2). Now noting that

$$K_1 - \Lambda K = G, \qquad \frac{d}{d\gamma}(K_1 - \Lambda K) = 0,$$

$$NK = \Sigma^*(\Sigma\Sigma^*)^{-1}A, \qquad \frac{d}{d\gamma}(NK) = 0,$$

$$K^*NK + C^*C = A^*(\Sigma\Sigma^*)^{-1}A \quad \text{and} \quad \frac{d}{d\gamma}(K^*NK + C^*C) = 0.$$

Then, by differentiating (4.14) with respect to $\gamma$, we obtain

$$(4.18) \quad \begin{aligned} &(K_1 + \Lambda N^{-1} \Lambda^* \overline{P})^* \frac{d\overline{P}}{d\gamma} + \frac{d\overline{P}}{d\gamma}(K_1 + \Lambda N^{-1} \Lambda^* \overline{P}) \\ &+ \frac{1}{(1-\gamma)^2}(\Lambda^* \overline{P} + NK)^* \Sigma^*(\Sigma\Sigma^*)^{-1}\Sigma(\Lambda^* \overline{P} + NK) = 0. \end{aligned}$$

Set $(\frac{d\overline{P}}{d\gamma})_{(-\infty)} = \lim_{\gamma \to -\infty} \frac{d\overline{P}}{d\gamma}$. Then, sending $\gamma$ to $-\infty$ in (4.18), we see that

$$(G + \Lambda \widehat{N}_{(-\infty)} \Lambda^* \overline{P}_{(-\infty)})^* \left( \frac{d\overline{P}}{d\gamma} \right)_{(-\infty)}$$

$$+ \left( \frac{d\overline{P}}{d\gamma} \right)_{(-\infty)} (G + \Lambda \widehat{N}_{(-\infty)} \Lambda^* \overline{P}_{(-\infty)}) = 0.$$

Since $G + \Lambda \widehat{N}_{(-\infty)} \Lambda^* \overline{P}_{(-\infty)}$ is stable, we see that

$$(4.19) \qquad \left( \frac{d\overline{P}}{d\gamma} \right)_{(-\infty)} = 0.$$

Set $\overline{q}_{(-\infty)} = \lim_{\gamma \to -\infty} \overline{q}(\gamma)$. As for $\overline{q}(\gamma)$, sending $\gamma$ to $-\infty$ in (2.35), we have

$$(G + \Lambda \widehat{N}_{(-\infty)} \Lambda^* \overline{P}_{(-\infty)})^* \overline{q}_{(-\infty)} + b^* \overline{P}_{(-\infty)}$$

$$- (A^* + \overline{P}_{(-\infty)} \Lambda \Sigma^*)(\Sigma\Sigma^*)^{-1}(a - r\mathbf{1}) = 0,$$

and so

$$\overline{q}_{(-\infty)} = \{(G + \Lambda \widehat{N}_{(\infty)} \Lambda^* \overline{P}_{(-\infty)})^*\}^{-1}$$
$$\times [(A^* + \overline{P}_{(-\infty)} \Lambda \Sigma^*)(\Sigma \Sigma^*)^{-1}(a - r\mathbf{1}) - b^* \overline{P}_{(-\infty)}].$$

Moreover, setting $(\frac{d\overline{q}}{d\gamma})_{(-\infty)} = \lim_{\gamma \to -\infty} \frac{d\overline{q}}{d\gamma}$, we see by (4.12) that

$$(G + \Lambda \widehat{N}_{(-\infty)} \Lambda^* \overline{P}_{(-\infty)})^* \left(\frac{d\overline{q}}{d\gamma}\right)_{(-\infty)} = 0.$$

Since $G + \Lambda \widehat{N}_{(-\infty)} \Lambda^* \overline{P}_{(-\infty)}$ is stable, we have

(4.20)
$$\left(\frac{d\overline{q}}{d\gamma}\right)_{(-\infty)} = 0.$$

The present proposition is directly seen by using (4.13), (4.19) and (4.20).
$\square$

**Acknowledgments.**    The authors would like to thank the referees for helpful comments and suggestions. The first author is a postdoctoral fellow at Institute of Mathematics, Academia Sinica, Taiwan. Therefore, he is thankful to Academia Sinica for constant support.

## REFERENCES

[1] BUCY, R. S. and JOSEPH, P. D. (1987). *Filtering for Stochastic Processes with Applications to Guidance*, 2nd ed. Chelsea Publishing, New York. MR0879420

[2] BIELECKI, T. R. and PLISKA, S. R. (1999). Risk-sensitive dynamic asset management. *Appl. Math. Optim.* **39** 337–360. MR1675114

[3] BIELECKI, T. R. and PLISKA, S. R. (2004). Risk-sensitive ICAPM with application to fixed-income management. *IEEE Trans. Automat. Control* **49** 420–432. MR2062254

[4] BIELECKI, T. R., PLISKA, S. R. and SHEU, S.-J. (2005). Risk sensitive portfolio management with Cox–Ingersoll–Ross interest rates: The HJB equation. *SIAM J. Control Optim.* **44** 1811–1843. MR2193508

[5] BROWNE, S. (1999). Beating a moving target: Optimal portfolio strategies for outperforming a stochastic benchmark. *Finance Stoch.* **3** 275–294. MR1842287

[6] BROWNE, S. (1999). The risk and rewards of minimizing shortfall probability. *J. Portfolio Management* **25** 76–85.

[7] DAVIS, M. and LLEO, S. (2008). Risk-sensitive benchmarked asset management. *Quant. Finance* **8** 415–426. MR2435642

[8] DEMBO, A. and ZEITOUNI, O. (1998). *Large Deviations Techniques and Applications*, 2nd ed. *Applications of Mathematics* **38**. Springer, New York. MR1619036

[9] EKELAND, I. and TÉMAM, R. (1999). *Convex Analysis and Variational Problems*, English ed. *Classics in Applied Mathematics* **28**. SIAM, Philadelphia, PA. MR1727362

[10] FLEMING, W. H. (1995). Optimal investment models and risk-sensitive stochastic control. In *Mathematical Finance* (M. DAVIS ET AL., eds.) 75–88. Springer, Berlin.

[11] FLEMING, W. H. and JAMES, M. R. (1995). The risk-sensitive index and the $H_2$ and $H_\infty$ norms for nonlinear systems. *Math. Control Signals Systems* **8** 199–221. MR1387043

[12] FLEMING, W. H. and MCENEANEY, W. M. (1995). Risk-sensitive control on an infinite time horizon. *SIAM J. Control Optim.* **33** 1881–1915. MR1358100

[13] FLEMING, W. H. and RISHEL, R. W. (1975). *Deterministic and Stochastic Optimal Control.* Springer, Berlin. MR0454768

[14] FLEMING, W. H. and SHEU, S.-J. (1999). Optimal long term growth rate of expected utility of wealth. *Ann. Appl. Probab.* **9** 871–903. MR1722286

[15] FLEMING, W. H. and SHEU, S. J. (2000). Risk-sensitive control and an optimal investment model. *Math. Finance* **10** 197–213. MR1802598

[16] FLEMING, W. H. and SHEU, S. J. (2002). Risk-sensitive control and an optimal investment model. II. *Ann. Appl. Probab.* **12** 730–767. MR1910647

[17] FLEMING, W. H. and SONER, H. M. (2006). *Controlled Markov Processes and Viscosity Solutions*, 2nd ed. *Stochastic Modelling and Applied Probability* **25**. Springer, New York. MR2179357

[18] HATA, H. and IIDA, Y. (2006). A risk-sensitive stochastic control approach to an optimal investment problem with partial information. *Finance Stoch.* **10** 395–426. MR2244352

[19] HATA, H. and SEKINE, J. (2005). Solving long term optimal investment problems with Cox–Ingersoll–Ross interest rates. *Adv. Math. Econ.* **8** 231–255.

[20] HÖRMANDER, L. (1967). Hypoelliptic second-order differential equations. *Acta Math.* **119** 147–171. MR0222474

[21] KAISE, H. and SHEU, S. J. (2004). Risk sensitive optimal investment: Solutions of the dynamical programming equation. In *Mathematics of Finance. Contemporary Mathematics* **351** 217–230. Amer. Math. Soc., Providence, RI. MR2076543

[22] KAISE, H. and SHEU, S.-J. (2006). On the structure of solutions of ergodic type Bellman equation related to risk-sensitive control. *Ann. Probab.* **34** 284–320. MR2206349

[23] KAISE, H. and SHEU, S.-J. (2006). Evaluation of large time expectations for diffusion processes. Preprint.

[24] KORN, R. (1997). *Optimal Portfolios.* World Scientific, Singapore.

[25] KURODA, K. and NAGAI, H. (2002). Risk-sensitive portfolio optimization on infinite time horizon. *Stoch. Stoch. Rep.* **73** 309–331. MR1932164

[26] MERTON, R. C. (1971). Optimum consumption and portfolio rules in a continuous-time model. *J. Econom. Theory* **3** 373–413. MR0456373

[27] NAGAI, H. (1996). Bellman equations of risk-sensitive control. *SIAM J. Control Optim.* **34** 74–101. MR1372906

[28] NAGAI, H. (2003). Optimal strategies for risk-sensitive portfolio optimization problems for general factor models. *SIAM J. Control Optim.* **41** 1779–1800. MR1972534

[29] NAGAI, H. and PENG, S. (2002). Risk-sensitive dynamic portfolio optimization with partial information on infinite time horizon. *Ann. Appl. Probab.* **12** 173–195. MR1890061

[30] PHAM, H. (2003). A large deviations approach to optimal long term investment. *Finance Stoch.* **7** 169–195. MR1968944

[31] PHAM, H. (2003). A risk-sensitive control dual approach to a large deviations control problem. *Systems Control Lett.* **49** 295–309. MR2011524

[32] SEKINE, J. (2007). *Mathematical Finance.* Baihukan, Japan.

[33] SEKINE, J. (2009). Private communication.

[34] RISHEL, R. (1999). Optimal portfolio management with partial observations and power utility function. In *Stochastic Analysis, Control, Optimization and Applications. Systems Control Found. Appl.* 605–619. Birkhäuser, Boston, MA. MR1702984

[35] STUTZER, M. (2000). A portfolio performance index. *Financial Analysts Journal* **56** 52–61.

[36] STUTZER, M. (2003). Portfolio choice with endogenous utility: A large deviations approach. *J. Econometrics* **116** 365–386. MR2002529

[37] WHITTLE, P. (1990). *Risk-Sensitive Optimal Control.* Wiley, Chichester. MR1093001

[38] WONHAM, W. M. (1968). On a matrix Riccati equation of stochastic control. *SIAM J. Control Optim.* **6** 681–697. MR0239161

[39] WONHAM, W. M. (1985). *Linear Multivariable Control: A Geometric Approach*, 3rd ed. *Applications of Mathematics* **10**. Springer, New York. MR0770574

H. HATA
S.-J. SHEU
ACADEMIA SINICA
NANKANG, TAIPEI 11529
TAIWAN
E-MAIL: hata@math.sinica.edu.tw
         sheusj@math.sinica.edu.tw

H. NAGAI
OSAKA UNIVERSITY
TOYONAKA
560-8531, OSAKA
JAPAN
E-MAIL: nagai@sigmath.es.osaka-u.ac.jp